\documentclass[11pt]{article}

\usepackage{IEEEtrantools}
\usepackage{latexsym}
\usepackage{amsfonts}
\usepackage{amssymb}
\usepackage{amsmath,amsfonts,amsthm}
\usepackage{mathrsfs}
\usepackage{multirow}
\usepackage{tikz}
 \usetikzlibrary{arrows,decorations.markings}

\newcommand{\bbbt}{\mathbb{T}}
\newcommand{\scrt}{\mathscr{T}}

\newcommand{\be}{\begin{equation}}
\newcommand{\ee}{\end{equation}}
\newcommand{\bea}{\begin{eqnarray}}
\newcommand{\eea}{\end{eqnarray}}
\newcommand{\bean}{\begin{eqnarray*}}
\newcommand{\eean}{\end{eqnarray*}}
\newcommand{\brray}{\begin{array}}
\newcommand{\erray}{\end{array}}
\newcommand{\biearray}{\begin{IEEEarray}{rCl}}
\newcommand{\eiearray}{\end{IEEEarray}}

\newcommand{\newsection}[1]{\setcounter{equation}{0}
\setcounter{dfn}{0}
\section{#1}}

\newtheorem{dfn}{Definition}[section]
\newtheorem{thm}[dfn]{Theorem}
\newtheorem{lmma}[dfn]{Lemma}
\newtheorem{ppsn}[dfn]{Proposition}
\newtheorem{crlre}[dfn]{Corollary}
\newtheorem{xmpl}[dfn]{Example}
\newtheorem{rmrk}[dfn]{Remark}

\newcommand{\bdfn}{\begin{dfn}\rm}
\newcommand{\bthm}{\begin{thm}}
\newcommand{\blmma}{\begin{lmma}}
\newcommand{\bppsn}{\begin{ppsn}}
\newcommand{\bcrlre}{\begin{crlre}}
\newcommand{\bxmpl}{\begin{xmpl}}
\newcommand{\brmrk}{\begin{rmrk}\rm}

\newcommand{\edfn}{\end{dfn}}
\newcommand{\ethm}{\end{thm}}
\newcommand{\elmma}{\end{lmma}}
\newcommand{\eppsn}{\end{ppsn}}
\newcommand{\ecrlre}{\end{crlre}}
\newcommand{\exmpl}{\end{xmpl}}
\newcommand{\ermrk}{\end{rmrk}}

\newcommand{\bbc}{\mathbb{C}}
\newcommand{\bbz}{\mathbb{Z}}
\newcommand{\bbn}{\mathbb{N}}
\newcommand{\bbr}{\mathbb{R}}

\newcommand{\cla}{\mathcal{A}}
\newcommand{\clb}{\mathcal{B}}
\newcommand{\clh}{\mathcal{H}}

\newcommand{\clk}{\mathcal{K}}
\newcommand{\cll}{\mathcal{L}}

\newcommand{\prf}{\noindent{\it Proof\/}: }

\def \qed { \mbox{}\hfill
$\Box$\vspace{1ex}}

\begin{document}

\author{{\sc Bipul Saurabh}}
\title{ On quantum quaternion spheres}
\maketitle
 \begin{abstract}  
	For the quantum symplectic group $SP_q(2n)$, we describe the $C^*$-algebra of continuous functions on the
	 quotient space $SP_q(2n)/SP_q(2n-2)$ as a universal $C^*$-algebra given by a finite set of generators and relations.
	 The proof involves a careful analysis of the relations, and use of the branching rules for representations of the 
	 symplectic group due to Zhelobenko. We then exhibit a set of generators of the $K$-groups of this $C^*$-algebra 
	 in terms of  generators of the $C^*$-algebra. 
\end{abstract}
{\bf AMS Subject Classification No.:} {\large 58}B{\large 34}, {\large
46}L{\large 87}, {\large
  19}K{\large 33}

\newsection{Introduction}
  The theory of quantum groups was first studied in the topological setting independently by Woronowicz~\cite{Wor-1987aa} and
  Vaksman \& Soibelman~\cite{VakSou-1988aa}. Both of these papers studied 
 the quantum $SU(2)$ group. Later Woronowicz developed the theory of compact quantum groups and their representation theory. 
 The notions of quantum subgroups and quantum homogeneous spaces were soon introduced by Podle\'s~(\cite{Pod-1995aa}).
 The most well-known example of compact quantum group is the $q$-deformation of the $SU(n)$ group whose representation theory 
 was obtained by Vaksman \& Soibelman (see~\cite{VakSou-1990ab}). Later, Korogodski \& Soibelman (\cite{KorSoi-1998aa})
 gave a complete classification of all the irreducible representations of the $C^*$-algebra $C(G_q)$ 
 where $G_q$ is the $q$-deformation of a classical simply connected semisimple compact Lie group.
 Vaksman \& Soibelman ~\cite{VakSou-1990ab} further studied quantum homogeneous space $SU_q(n)/SU_q(n-1)$
 and obtained its irreducible representations. 
 Quantum analogs of the Stiefel manifolds $SU_{q}(n)/SU_q(n-m)$ were introduced by Podkolzin \& Vainerman~\cite{PodVai-1999aa}
 who classified  the irreducible representations of the $C^{*}$-algebras underlying such manifolds. 

Given a compact quantum group $G$ and a subgroup $H$, the $C^*$-algebra $C(G/H)$ of the quotient space $G/H$
is defined to be a certain $C^*$-subalgebra of $C(G)$. 
Extending the ideas of Korogodski \& Soibelman in (\cite{KorSoi-1998aa}), Neshveyev \& Tuset (\cite{NesTus-2012ab}) gave a 
complete classification of the irreducible representations of the $C^*$-algebra $C(G_q/H_q)$ where
$G_q$ is the $q$-deformation of a simply connected semisimple compact Lie group and $H_q$ is
the $q$-deformation of a closed  Poisson-Lie subgroup  $H$ of $G$.
To understand the structure of these $C^{*}$-algebras further, the next thing to look for would be their $K$-groups. 
Vaksman \& Soibelman~\cite{VakSou-1990ab}  determined $K$-groups of the stiefel manifolds $SU_{q}(n)/SU_{q}(n-1)$ which 
are quantum analogues of the classical odd dimensional spheres. 
Chakraborty \& Sundar\cite{ChaSun-2011aa} computed the $K$-theory of the Stiefel manifolds $SU_{q}(n)/SU_{q}(n-2)$ and 
as a consequence they found the $K$-theory of $SU_{q}(3)$.
Then Neshveyev \& Tuset proved in the above mentioned paper (\cite{NesTus-2012ab}) that
$C(G_q/H_q)$ is $KK$-equivalent to the classical counterpart $C(G/H)$.

Here one question arises. Can we describe the algebra of functions 
on a quantum group $G_q$ or its quotient spaces $G_q/H_q$ in terms of a finite set of generators and relations? Given a semisimple compact 
Lie group $G$, the algebra of functions $C(G_q)$ on its  $q$-deformation $G_q$ can be defined as the  universal
enveloping $C^*$-algebra generated by 
matrix co-efficients of all finite dimensional representations of Quantized universal enveloping algebra $U_q(g)$
where $g$ is the Lie algebra of $G$. 
 Similarly, for a subgroup $H_q$ of $G_q$, the algebra of functions $C(G_q/H_q)$ on the quotient space $G_q/H_q$
can be described as the subalgebra of $C(G_q)$ generated by matrix elements of certain finite dimensional representations
of Quantized universal enveloping algebra $U_q(g)$. But  how to get a finite set of generators and relations to describe these 
$C^*$- algebras is not clear. Korogodski and Soibelman described  $C(SU_q(2))$ in terms of finite set of 
generators and relations (see \cite{KorSoi-1998aa} theorem 3.0.1). For $C(SU_q(n))$ case, it is proved by  Woronowicz~\cite{Wor-1988aa}. 
Another approach called FRT approach to define quantum deformation of a simple Lie group is 
due to Faddeev, Reshetikhin and Takhtajan. They started with a finite set of generators and relations based on quantum $R$-matrix of
 a simple Lie algebra to define the algebra of functions on the quantum group. Through a pairing between this algebra of functions and 
 Quantized universal enveloping algebra, one can view this algebra of functions as certain subalgebra of dual Hopf-algebra of 
 Quantized universal enveloping algebra. This relates FRT approach with the previous one due to V. Drinfeld. One difficulty with FRT approach is 
  that there is no formula for the quantum $R$-matrix of an arbitrary simple Lie algebra. Also, even in this approach, the question of describing 
 the algebra of functions on homogeneous spaces does not have an answer yet.

Many of the $C^*$-algebras arising in the context of quantum groups and their homogeneous spaces
are universal $C^*$-algebras given by finite sets of generators and relations. This fact along with a knowledge of their generators and
relations turn out to be extremely useful in studying these objects further as well as
in doing various computations involving them.
For example, the description of $C(SU_q(2))$ in terms of its generators and relations was used by Vaksman \& Soibelman~\cite{VakSou-1990ab} 
to compute its K-theory. Later these generators and relations were used by various authors (\cite{ChaPal-2003aa}, \cite{Con-2004aa},
\cite{DabLanSit-2005aa}, \cite{PalSun-2010aa}) to study the quantum group $SU_q(2)$ in the context of non commutative geometry.
Similarly the description of the quantum odd dimensional spheres as universal $C^*$-algebras
given by finite sets of generators and relations was used by 
Vaksman \& Soibelman~\cite{VakSou-1990ab} to compute its K-theory; subsequently by Hong \& Szymanski~\cite{HonSzy-2002aa} 
to show that it can be realized as a graph $C^*$-algebra, which in turn gives a great deal of insight into these $C^*$-algebras,
and finally by Chakraborty \& Pal \cite{ChaPal-2008aa} and by 
Pal \& Sundar \cite{PalSun-2010aa} to study spectral triples for these spaces. 

Describing a $C^*$-algebra in terms of generators and relations is a rather difficult
and intractable problem in general. For $C^*$-algebra of 'continuous functions' on quantum groups or 
their homogeneous spaces, it is often the case that they are 
generated by certain matrix elements of representations. 
But there is no general results in this direction. In the present paper, we attempt this problem for the $C^*$-algebra  corresponding to 
the quantum Stiefel manifold 
$SP_q(2n)/SP_q(2n-2)$. We first show that 
$C(SP_q(2n)/SP_q(2n-2))$ is the $C^*$-subalgebra of $C(SP_q(2n))$ generated by elements of the first and last row
of the fundamental co-representation of the quantum symplectic group $SP_q(2n)$. 
We then get a finite set of relations satisfied by these elements. Then making use of the results of Neshveyev \& Tuset, we prove that 
$C(SP_q(2n)/SP_q(2n-2))$ is the 
universal $C^*$-algebra given by a finite set of generators and relations.
We then go on to exhibit generators of its $K$-groups in terms of the generators of this $C^*$-algebra.

Here is brief outline of how the paper is organized.
In subsection~2.1, we start by recalling $q$-deformation of the symplectic group $SP_{q}(2n)$ at the 
$C^{*}$-algebra level and describe its quantum group structure.
In subsection $(2.2)$ and $(2.3)$, we first write down explicitly the pairing between
the quantum universal enveloping algebra $U_q(sp(2n))$ and the co-ordinate algebra
$O(SP_q(2n))$, which helps us write down all elementary representations. An application of the result of Korogodski \& Soibelman then gives us
 a complete list of all irreducible representations of $C(SP_{q}(2n))$, which in turn
 provides us with a faithful representation on a Hilbert space. In subsection $2.4$, extending the ideas of Chakraborty \& Pal~\cite{ChaPal-2003ac},
 we associate 
 certain diagrams with these irreducible representations.

In section~3, we consider the quantum homogeneous space $SP_q(2n)/SP_q(2n-2)$ and 
show that the $C^*$-algebra underlying this manifolds is generated by a finite set of generators. Here we use branching rules for representations
of the symplectic group due to Zhelobenko~\cite{Zhe-1962aa}.
In section~4, we define the $C^*$-algebra $C(H_q^{2n})$ of continuous functions on the quantum quaternion sphere
as a universal $C^{*}$-algebra with $2n$ generators satisfying a certain set of  relations. 
We then study these relations and find all irreducible representations of the $C^*$-algebra $C(H_q^{2n})$. 
We then show that the $C^*$-algebra $C(SP_{q}(2n)/SP_{q}(2n-2))$ is a homomorphic image of $C(H_q^{2n})$. This is 
accomplished by showing that the generators of $C(SP_{q}(2n)/SP_{q}(2n-2))$ obtained in section~3 satisfy 
the defining relations of $C(H_q^{2n})$. Finally, by comparing the irreducible representations of the two $C^*$-algebras,
we conclude that they are isomorphic. Here we use Neshveyev \& Tuset's results on irreducible representations of a quotient space.
Finally, in section~5, we find a chain of short exact sequences for $C\left(SP_{q}(2n)/SP_{q}(2n-2)\right)$. Utilizing these,
we compute their $K$-groups with explicit generators. Here we should remark that these $K$-groups are already known, 
thanks to the $KK$-equivalence of $C\left(SP_{q}(2n)/SP_{q}(2n-2)\right)$ and  $C\left(SP(2n)/SP(2n-2)\right)$ proved 
in \cite{NesTus-2012ab}. But the advantage of our computation here is that we produce explicit generators of these $K$-groups 
which can come in handy in many situations.

A word about notations. It is unfortunate but very common to use the notation $sp(2n, \mathbb{C})$ to describe the type $C_n$ groups
 at the Lie algebra level, while one switches to $SP(n)$ at the group level. Here we have used $2n$ at both places, i.e. $U_q(sp(2n))$
 denotes the Quantized universal enveloping algebra of type $C_n$ and $SP_q(2n)$ denotes the same quantum group at the function algebra 
 level. Throughout this paper, $q$ will denote a real number in the interval  $(0,1)$ and $C$ is used to denote a generic constant. 
 The standard bases of the 
Hilbert spaces $\ell^2(\bbn)$ and $\ell^2(\bbz)$ will be denoted by
$\left\{e_n: n\in \bbn \right\}$ and  $\left\{e_n: n\in \bbz \right\}$ respectively.

\newsection{Quantum symplectic group }

\subsection{The $C^*$-algebra $C(SP_q(2n))$} 
In this subsection, we briefly recall from \cite{KliSch-1997aa} various algebras and maps related to the compact
quantum group $SP_q(2n)$. We set up some notation that will be used throughout this paper. Define, 
\begin{IEEEeqnarray}{rCl} 
i^{'} & = & 2n+1-i \nonumber \\
\rho_{i}  & = & n+1-i \quad \text{ if } i \leq n. \nonumber \\
\rho_{i^{'}} &=& -\rho_{i}. \nonumber \\
\epsilon_{i}&=&\begin{cases}
               1 & \mbox{ if } 1\leq i\leq n,\cr
              -1 & \mbox{ if } n+1\leq i\leq 2n,\cr
               \end{cases} \nonumber \\         
C_{j}^{i} &=&\epsilon_{i}\delta_{ij} q^{-\rho_{i}}. \nonumber \\
\theta(i)&=&\begin{cases}
               0 & \mbox{ if }  i\leq 0,\cr
               1 & \mbox{ if }  i > 0,\cr
               \end{cases} \nonumber \\
R_{mn}^{ij} &=& q^{\delta_{ij}-\delta_{ij^{'}}}\delta_{im}\delta_{jn} + (q-q^{-1})\theta(i-m)(\delta_{jm}\delta_{in}+ C_{j}^{i}C_{n}^{m}). \nonumber
\end{IEEEeqnarray}
Let $\bbc\left\langle u_{j}^{i}\right\rangle$ denote the free algebra with generators $u_{j}^{i}$, $i,j=1,2,\cdots 2n$ and
let $J(R)$ be the two sided ideal of $\bbc\left\langle u_{j}^{i}\right\rangle$ generated by the following elements,
\begin{IEEEeqnarray}{rCl} 
                    I_{st}^{ij}=\sum_{k,l=1}^{2n}R_{kl}^{ji}u_{s}^{k}u_{t}^{l}-R_{st}^{lk}u_{k}^{i}u_{l}^{j},
                    \quad i,j,s,t=1,2\cdots 2n. \label{relations}
\end{IEEEeqnarray}
Let $A(R)$ denote the quotient algebra $\bbc\left\langle u_{j}^{i}\right\rangle/J(R)$. The $2n \times 2n$ matrices $(\!(u_{j}^{i})\!)$
and $(\!(C_{j}^{i})\!)$ are denoted by $U$ and $C$ respectively. Define $J =\left\langle UCU^{t}C^{-1}-I,CU^{t}C^{-1}U-I\right\rangle$, the
two sided ideal generated by entries of matrices $UCU^{t}C^{-1}-I$ and $CU^{t}C^{-1}U-I$. Let $\emph{O}(SP_{q}(2n))$ denote the quotient 
algebra $A(R)/J$.

The algebra $\emph{O}(SP_{q}(2n))$ is a Hopf-*algebra with co-multiplication $\Delta$, co-unit $\epsilon$, antipode $S$ and
involution $*$ given on the generating elements by,
\begin{displaymath}
\Delta(u_{j}^{i}) = \sum_{k=1}^{N}u_{k}^{i} \otimes u_{j}^{k},\qquad 
\epsilon(u_{j}^{i}) = \delta_{ij}, \qquad
S(u_{j}^{i}) = \epsilon_{i}\epsilon_{j}q^{\rho_{i}-\rho_{j}}u_{i^{'}}^{j^{'}}, \qquad
(u_{j}^{i})^{*} = \epsilon_{i}\epsilon_{j}q^{\rho_{i}-\rho_{j}}u_{j^{'}}^{i^{'}}.
\end{displaymath}
Note that $U^{*} = CU^{t}C^{-1}$. Hence we have,
\be \label{e}
UU^{*}=U^{*}U=I.
\ee
Now to make $\emph{O}(SP_{q}(2n))$, a normed-$*$algebra, we define,
\[
     \left\|a\right\| = \mbox{ sup }\left\{\left\|\pi(a)\right\|: \pi \mbox{ is a representation of }\emph{O}(SP_{q}(2n))\right\}.
\]
By $(\ref{e})$, we have, $\left\|u_{j}^{i}\right\|\leq 1$, hence for all $a\in \emph{O}(SP_{q}(2n))$,  $\left\|a\right\|<\infty$.
We denote by $C(SP_{q}(2n))$  the completion of $\emph{O}(SP_{q}(2n))$. The pair $(C(SP_{q}(2n)),\Delta)$ is a compact quantum group
called a $q$-deformation of the symplectic group $SP_{q}(2n)$. For more detail, see (\cite{KliSch-1997aa}, page~317--320, page~326).

\subsection{Pairing between $U_{q}(sp_{2n})$ and $\emph{O}(SP_{q}(2n))$}  
Let $(\!(a_{ij})\!)$ be the Cartan matrix of Lie algebra $sp(2n)$ given by,
\[
a_{ij}=\begin{cases}
               2 & \mbox{ if } i=j,\cr
              -1 & \mbox{ if } i=j+1,\cr
							-1 & \mbox{ if } i=j-1,i\neq n-1,\cr
							-2 & \mbox{ if } i=j-1=n-1,\cr
               0 & otherwise,
               \end{cases}
\]
Define $q_{i}=q^{d_{i}}$, where $d_{i}$= 1 for $i=1,2,..n-1$ and $d_{n}=2$. 
The quantized universal enveloping algebra (QUEA) $U_{q}(sp_{2n})$ is the universal algebra
generated by $E_i$, $F_i$, $K_i$ and $K_i^{-1}$, $i=1,\ldots,\ell$, satisfying the
following relations
\begin{displaymath}
K_iK_j=K_jK_i,\quad K_iK_i^{-1}=K_i^{-1}K_i=1,
\end{displaymath}
\begin{displaymath}
K_iE_jK_i^{-1}=q_i^{a_{ij}}E_j,\quad
K_iF_jK_i^{-1}=q_i^{ a_{ij}}F_j,
\end{displaymath}
\begin{displaymath}
E_iF_j-F_jE_i=\delta_{ij}\frac{K_i-K_i^{-1}}{q_i-q_i^{-1}},
\end{displaymath}
\begin{displaymath}
\sum_{r=0}^{1-a_{ij}}(-1)^r 
     \binom{1-a_{ij}}{r}_{q_i}
  E_i^{1-a_{ij}-r}E_jE_i^r =0 \quad\forall\, i\neq j,
\end{displaymath}
\begin{displaymath}
\sum_{r=0}^{1-a_{ij}}(-1)^r 
     \binom{1-a_{ij}}{r}_{q_i}
  F_i^{1-a_{ij}-r}F_jF_i^r =0\quad \forall\, i\neq j,
\end{displaymath}
where $\binom{n}{r}_q$ denote the $q$-binomial coefficients.
Hopf *-structure comes from the following maps:
\[
\Delta(K_i)=K_i\otimes K_i,\quad \Delta(K_i^{-1})=K_i^{-1}\otimes K_i^{-1},
\]
\[
\Delta(E_i)=E_i\otimes K_i  + 1\otimes E_i,\quad
\Delta(F_i)=F_i\otimes 1  + K_i^{-1}\otimes F_i,
\]
\[
\epsilon(K_i)=1,\quad \epsilon(E_i)=0=\epsilon(F_i),
\]
\[
S(K_i)=K_i^{-1},\quad S(E_i)=-E_iK_i^{-1},\quad S(F_i)=-K_iF_i,
\]
\[
K_i^*=K_i,\quad E_i^*=K_iF_i,\quad F_i^*=E_iK_i^{-1}.
\]
We refer the reader to \cite{KliSch-1997aa} for a proof of the following theorem that gives a dual pairing between
$U_{q}(sp_{2n})$ and $\emph{O}(SP_{q}(2n))$.
\bthm $(\cite{KliSch-1997aa})$ \label{pairing}
There exist unique dual pairing $\left\langle\cdot,\cdot\right\rangle$ between the Hopf algebras and $U_{q}(sl_{2})$ 
and $\emph{O}(SL_{q}(2))$ and between $U_{q}(sp_{2n})$ and $\emph{O}(SP_{q}(2n))$ such that
\[
  \left\langle f,u_{l}^{k}\right\rangle = t_{kl}(f)        \hspace{1in} \mbox{for }  k,l=1,2,...2n.\\
\]
where $t_{kl}$ is the matrix element of $T_{1}$, the vector representation of $U_{q}(sl_{2})$ in first case and 
that of  $U_{q}(sp_{2n})$ in second case.
\ethm

We will explicitly describe $T_{1}$ for both cases and determine the pairing. For that, let $E_{ij}$ be the $2n \times 2n$ 
matrix with $1$ in the $(i,j)^{th}$ position and $0$ elsewhere and $D_{j}$ be the diagonal matrix with $q$ in the $(j,j)^{th}$ 
position and $1$ elsewhere on the diagonal.\\
For the QUEA $U_{q}(sp_{2n})$, one has
\[ \left.
\begin{array}{rcl}
T_1(K_i) &=& D_i^{-1}D_{i+1}D_{2n-i}^{-1}D_{2n-i+1}. \nonumber\\
T_1(E_i) &=& E_{i+1,i}-E_{2n-i+1,2n-i}.\nonumber\\
T_1(F_i) &=& E_{i,i+1}-E_{2n-i,2n-i+1}.\nonumber
\end{array} \right\}\quad  \mbox{for} \quad i \in \{1,2,\cdots,n-1\}
\]
and 
\[ \left.
\begin{array}{rcl}
T_1(K_n) = D_n^{-2}D_{n+1}^{2}, \qquad
T_1(E_n) = E_{n+1,n}, \qquad
T_1(F_n) = E_{n,n+1} \nonumber
\end{array} \right\} \quad \mbox{for}\quad i=n.
\]

For $U_{q}(sl_{2})$, one has (here $E_{ij}$ and $D_{j}$ are the $2 \times 2$ matrices.) 
\begin{IEEEeqnarray}{rCl}
T_1(K) = D_1^{-1}D_{2},\qquad
T_1(E) = E_{2,1},\qquad
T_1(F) = E_{1,2}.\nonumber
\end{IEEEeqnarray}
 For more detail, we refer to (\cite{KliSch-1997aa}, page~267--268). We will use these two pairing to 
 write down irreducible representations of $\emph{O}(SP_{q}(2n))$ which can be extended to $C(SP_{q}(2n))$ to get
 elementary representations of $C(SP_{q}(2n))$. \\
\textbf{Elementary representation of $C(SP_{q}(2n))$:} For $i \in \{1,2,\cdots,n\}$, let $\varphi_i : U_{q_i}(sl_2) \longrightarrow  U_q(sp(2n))$
be a  $*$-homomorphism given on generators of $U_{q_i}(sl_2)$ by, 
\begin{displaymath}
 K \longmapsto K_i, \qquad E \longmapsto E_i, \qquad F \longmapsto F_i.
\end{displaymath}
Consider the dual epimorphism, 
\begin{displaymath}
\varphi_i^* : C(SP_q(2n)) \longrightarrow C(SU_{q_i}(2))
\end{displaymath}
such that 
\begin{displaymath}
\left\langle f, \varphi_i^*(u_n^m)\right\rangle = \left\langle \varphi_i(f), u_n^m\right\rangle
\end{displaymath}
In particular,
\begin{IEEEeqnarray}{rCl}
\left\langle K, \varphi_i^*(u_n^m)\right\rangle = \left\langle K_i, u_n^m\right\rangle, \quad
\left\langle E, \varphi_i^*(u_n^m)\right\rangle = \left\langle E_i, u_n^m\right\rangle, \quad 
\left\langle F, \varphi_i^*(u_n^m)\right\rangle = \left\langle F_i, u_n^m\right\rangle. \label{r}
\end{IEEEeqnarray}
\brmrk Initially $\varphi_{i}$ will induce an $*$ epimorphism from the Hopf-$*$algebra $\emph{O}(SP_{q}(2n))$ to $\emph{O}(sl_{q_{i}}(2))$ 
which when extended to $C(SP_{q}(2n))$ gives the above homomorphism at $C^{*}$-algebra level.
For more detail see (\cite{KliSch-1997aa}, page~327) and \cite{KorSoi-1998aa}.
\ermrk 
Let $N$ be the number operator given by $N: e_{n} \mapsto n e_{n}$ 
and $S$ be the  shift operator given by $S: e_{n} \mapsto e_{n-1}$ on $L_{2}(\bbn)$. Denote by $\pi$ the following 
representation of $C(SU_{q}(2))$ on $L_{2}(\bbn)$;\\
\[
\pi(u_l^k)=\begin{cases}
              \sqrt{1-q^{2N+2}}S & \mbox{ if } k=l=1,\cr
              S^{*}\sqrt{1-q^{2N+2}} & \mbox{ if } k=l=2,\cr
							-q^{N+1} & \mbox{ if } k=1,l=2,\cr
							q^N & \mbox{ if } k=2,l=1,\cr
							\delta_{kl} & \mbox{ otherwise }. \cr
              \end{cases}
\] 
Define, $\pi_{s_{i}} = \pi \circ \varphi_{i}^{*}.$ Applying $(\ref{r})$, we have, for $i=1,2,\cdots,n-1$,
\[
\pi_{s_{i}}(u_l^k)=\begin{cases}
              \sqrt{1-q^{2N+2}}S & \mbox{ if } (k,l)=(i,i) \mbox{ or } (2n-i,2n-i),\cr
              S^*\sqrt{1-q^{2N+2}} & \mbox{ if } (k,l)=(i+1,i+1) \mbox{ or } (2n-i+1,2n-i+1),\cr
							-q^{N+1} & \mbox{ if } (k,l)=(i,i+1),\cr
							q^N & \mbox{ if } (k,l)=(i+1,i),\cr
							q^{N+1} & \mbox{ if } (k,l)=(2n-i,2n-i+1),\cr
							-q^N & \mbox{ if } (k,l)=(2n-i+1,2n-i),\cr
							\delta_{kl} & \mbox{ otherwise }. \cr
							\end{cases}
\]
For $i=n$,
\[
\pi_{s_n}(u_l^k)=\begin{cases}
              \sqrt{1-q^{4N+4}}S & \mbox{ if } (k,l)=(n,n),\cr
              S^*\sqrt{1-q^{4N+4}} & \mbox{ if } (k,l)=(n+1,n+1),\cr
							-q^{2N+2} & \mbox{ if } (k,l)=(n,n+1),\cr
							q^{2N} & \mbox{ if } (k,l)=(n+1,n),\cr
							\delta_{kl} & \mbox{ otherwise }. \cr
              \end{cases}
\]

Each $\pi_{s_{i}}$ is an  irreducible representation and is called an elementary representation of $C(SP_{q}(2n))$.
For any two representations $\varphi$ and $\psi$ of $C(SP_{q}(2n))$ define, $\varphi * \psi := (\varphi \otimes \psi)\circ \Delta$.
Let W be the Weyl group of $sp_{2n}$ and $\vartheta \in W$ such that  $s_{i_{1}}s_{i_{2}}...s_{i_{k}}$ is a reduced expression for $\vartheta$. 
Then $\pi_{\vartheta}= \pi_{s_{i_{1}}}*\pi_{s_{i_{2}}}*\cdots *\pi_{s_{i_{k}}}$ is an irreducible representation which is independent 
of the reduced expression. Now for $t=(t_{1},t_{2},\cdots ,t_{n}) \in \bbbt^{n}$, define the map  $\tau_t: C(SP_{q}(2n) \longrightarrow \bbc $ by
\[
\tau_{t}(u_j^i)=\begin{cases}
              \overline{t_{i}}\delta_{ij} & \mbox{ if } i \leq n,\cr
							t_{2n+1-i}\delta_{ij} & \mbox{ if } i > n,\cr
							\end{cases}
\]
Then $\tau_{t}$ is a $*$-algebra homomorphism. For $t \in \bbbt^{n}, \vartheta \in W$,  let $\pi_{t,\vartheta} = \tau_{t}*\pi_{\vartheta}$.

We refer to (\cite{KorSoi-1998aa}, page~121) for the following theorem.
\bthm
$\left\{\pi_{t,\vartheta}; t \in \bbbt^{n}, \vartheta \in W\right\}$ is a complete set of mutually inequivalent representations of $C(SP_{q}(2n))$.
\ethm

\subsection{Representations of $C(SP_{q}(2n))$}
To write down all irreducible representations of $C(SP_{q}(2n))$, we need to recall a few facts on  the Weyl group of $sp_{2n}$ which
we summarize below. Weyl group $W_{n}$ of $sp_{2n}$ is a Coxeter group  generated by $s_{1},s_{2},...s_{n}$ satisfying the following relations:
\begin{align*}
s_{i}^{2} &= 1 & \text{for } & i=1,2,...n\\
s_{i}s_{i+1}s_{i} &=s_{i+1}s_{i}s_{i+1} & \text{for } & i=1,2,...n-1\\
s_{n-1}s_{n}s_{n-1}s_{n} &= s_{n}s_{n-1}s_{n}s_{n-1}
\end{align*}
The group  $W_{n}$ can be embedded faithfully in $M_{n}(\bbr)$ as,
\begin{IEEEeqnarray}{rCll}
s_{i} &=& I-E_{i,i}-E_{i+1,i+1}+E_{i,i+1}+E_{i+1,i}, & \qquad \mbox{for}\quad i=1,2,...n-1,\nonumber \\
s_{n} &=& I-2E_{n,n}, & \qquad \mbox{for} \quad i=n.\nonumber
\end{IEEEeqnarray}
So, $W_{n}$ is isomorphic to a subgroup of $GL(n,\bbr)$ generated by $s_{1},s_{2},...s_{n}$. We refer to \cite{FulHar-1991aa} for a proof 
of the following proposition.

\bppsn \label{pw1}
\begin{enumerate}
\item
Let $\mathfrak{S}_{n}$ be the permutation group and $H_n$ be the $n$-fold direct product of the group 
$\left\{-1,1\right\}$. $\mathfrak{S}_{n}$ acts on $H_n$ by permuting its co-ordinates. Then $W_{n} = H_n\rtimes \mathfrak{S}_{n}$.\\
In other word, $W_{n}$ is the set of $n \times n$ matrices having one non-zero entry in each row and each column which is either $1$ or $-1$.
\item
Any element of $W_{n}$ can be written in the form: $\quad \prod_{r=1}^{n}\psi_{r,k_{r}}^{(\epsilon_{r})}$ where 
$\epsilon_{r} \in \left\{0,1,2\right\}$ and $r\leq k_{r}\leq n$ with the convention that,
\[
\psi_{r,k_{r}}^{\epsilon}=\begin{cases}
               s_{k_{r}-1}s_{k_{r}-2}...s_{r} & \mbox{ if } \epsilon=1,\cr
               s_{k_{r}}s_{k_{r}+1}\cdots...s_{n-1}s_{n}s_{n-1}\cdots s_{k_{r}}s_{k_{r}-1}\cdots s_{r} & \mbox{ if } \epsilon=2,\cr
							 \mbox{ empty string } & \mbox{ if } \epsilon=0,\cr
               \end{cases}
\]
Also, the above expression is a reduced expression.
\item
The longest word of $W_{n}$ is $-I$ which can be written as $\quad \prod_{r=1}^{n}\psi_{r,r}^{(2)}$.
Also, $\left\{\psi_{r,r}^{(2)}\right\}_{i=1}^{n}$ commutes, hence $-I$ can be written as $\quad \prod_{r=1}^{n}\psi_{n+1-r,n+1-r}^{(2)}$,
which is a reduced expression.\\
\end{enumerate}
\eppsn

Denote by $\scrt$ the Toeplitz algebra. Let $\vartheta$ be a word on $s_{1},s_{2},...s_{n}$ of length $\ell(\vartheta)$. 
Then the map $\bbbt^{n} \ni t\longmapsto \pi_{t,\vartheta}(u_{j}^{i})\in\scrt^{\otimes \ell(\vartheta)}$ is continuous.
Hence we have a homomorphism $\chi_{\vartheta} : C(SP_{q}(2n))\longrightarrow C(\bbbt^{n})\otimes \scrt^{\otimes \ell(\vartheta)}$ 
such that $\chi_{\vartheta}(a)(t) = \pi_{t,\vartheta}(a),  \mbox{ for all } a\in C(SP_{q}(2n))$.

\bppsn
If $\vartheta^{'}$ is a subword of $\vartheta$ then $\chi_{\vartheta_{'}}$ and $\pi_{t,\vartheta_{'}}$ factor through $\chi_{\vartheta}$.
\eppsn
The proof of above proposition is straightforward, so we omit it. 
\bthm
Let $\vartheta_{n}$ be the longest word of the Weyl group of $sp_{2n}$ i.e 
\[
 \vartheta_{n} = (s_{n})(s_{n-1}s_{n}s_{n-1})...(s_{2}...s_{n}..s_{2})(s_{1}s_{2}...s_{n-1}s_{n}s_{n-1}...s_{1}).
\]
  Then the homomorphism
\[
                            \chi_{\vartheta_{n}} : C(SP_{q}(2n))\longrightarrow C(\bbbt^{n})\otimes \scrt^{\otimes \ell(\vartheta_{n})}
\]
is faithful.
\ethm
\prf Any  irreducible representation of $C(SP_{q}(2n))$ is of the form $\pi_{t,\vartheta}$ where $\vartheta$ is a 
word on  $s_{1},s_{2},...s_{n}$ and $t \in \bbbt^{n}$. From proposition $\ref{pw1}$, it is clear that  $\vartheta$ is 
a subword of $\vartheta_{n}$, hence $\pi_{t,\vartheta}$ factors through $\chi_{\vartheta_{n}}$  which shows that
$\chi_{\vartheta_{n}}$ is faithful. \qed

\subsection{Diagram representation}

At this point it will be useful to associate some diagrams with the 
above representations. We will use the scheme followed by 
Chakraborty \& Pal~\cite{ChaPal-2003ac} with a few additions. 
For convenience, we use labeled arrows to 
represent operators as given in the following table.\\

         \begin{tabular}{|c|c||c|c||c|c|}
   \hline
   Arrow type & Operator &  Arrow type  & Operator &  Arrow type  & Operator \\
   \hline   
    &&
   \multirow{2}{*}{}&
  \multirow{2}{*}{} &
 \multirow{2}{*}{} &
    \multirow{2}{*}{}\\
\begin{tikzpicture}[scale=.7]
\draw [->] (0,0) -- (1,0);
\end{tikzpicture}
  &
  $I$ &
     &&&\\
     &&&&&\\
     \hline  
     &&&&&\\
 \begin{tikzpicture}[scale=.7]
 \draw [->] (0,0) -- (1,0);
 \node at (.7,.3) {$+$};
 \end{tikzpicture}
   &
   $S^*\sqrt{I-q^{2N+2}}$ &
\begin{tikzpicture}[scale=.7]
\draw [->] (0,0) -- (1,1);
\node at (.9,.5) {$+$};
\end{tikzpicture}
&
$-q^{N+1}$&
\begin{tikzpicture}[scale=.7]
\draw [->] (0,1) -- (1,0);
\node at (.9,.5) {$+$};
\end{tikzpicture}
       &
$q^{N}$\\
  &&&&&\\
     \hline
     &&&&&\\
\begin{tikzpicture}[scale=.7]
\draw [->] (0,0) -- (1,0);
\node at (.7,.3) {$-$};
\end{tikzpicture}
  &
  $\sqrt{I-q^{2N+2}}S$ &
\begin{tikzpicture}[scale=.7]
\draw [->] (0,0) -- (1,1);
\node at (.9,.5) {$-$};
\end{tikzpicture}
   &
$q^{N+1}$&
\begin{tikzpicture}[scale=.7]
\draw [->] (0,1) -- (1,0);
\node at (.9,.5) {$-$};
\end{tikzpicture}
	  &
$-q^{N}$\\
 &&&&&\\
     \hline
     &&&&&\\
\begin{tikzpicture}[scale=.7]
\draw [->] (0,0) -- (1,0);
\node at (.7,.3) {$++$};
\end{tikzpicture}
  &
  $S^*\sqrt{I-q^{4N+4}}$ &
\begin{tikzpicture}[scale=.7]
\draw [->] (0,0) -- (1,1);
\node at (1,.5) {$++$};
\end{tikzpicture}
		&
$-q^{2N+2}$&
\begin{tikzpicture}[scale=.7]
\draw [->] (0,1) -- (1,0);
\node at (1,.5) {$++$};
\end{tikzpicture}
	   &
 $q^{2N}$\\
  &&&&&\\
     \hline
     &&&&&\\
\begin{tikzpicture}[scale=.7]
\draw [->] (0,0) -- (1,0);
\node at (.7,.3) {$-\,-$};
\end{tikzpicture}
  &
  $\sqrt{I-q^{4N+4}}S$ &
     &&&\\
     &&&&&\\
     \hline
\end{tabular} \\

Let us describe how to use a diagram to represent the irreducible $\pi_{s_i}$.

\begin{tabular}{p{200pt}p{200pt}}
\begin{tikzpicture}[scale=1.2]
\draw [->] (0,5) -- (1,5);
\draw [dashed] (0,4.5) -- (1,4.5);
\draw [->] (0,4) -- (1,4);
\draw [->] (0,4) -- (1,3);
\draw [->] (0,3) -- (1,4);
\draw [->] (0,3) -- (1,3);
\draw [dashed] (0,2.5) -- (1,2.5);
\draw [->] (0,2) -- (1,2);
\draw [->] (0,2) -- (1,1);
\draw [->] (0,1) -- (1,2);
\draw [->] (0,1) -- (1,1);
\draw [dashed] (0,.5) -- (1,.5);
\draw [->] (0,0) -- (1,0);
\node at (.5,2.2) {${}+$};
\node at (1,1.7) {${}+$};
\node at (1,1.3) {${}+$};
\node at (1,3.7) {${}-$};
\node at (1,3.3) {${}-$};
\node at (.5,.85) {${}-{}$};
\node at (.5,4.2) {${}+$};
\node at (.5,2.8) {${}-$};
\node at (-.5,5.05){$2n$};
\node at (1.5,5.05){$2n$};
\node at (-.8,4.05){$2n-i+1$};
\node at (1.9,4.05){$2n-i+1$};
\node at (-.5,3.05){$2n-i$};
\node at (1.5,3.05){$2n-i$};
\node at (-.5,2.05){$i+1$};
\node at (1.5,2.05){$i+1$};
\node at (-.3,1.05){$i$};
\node at (1.3,1.05){$i$};
\node at (-.3,.05){$1$};
\node at (1.3,.05){$1$};
 \node at (.5,-1){\text Diagram 1: $\psi_{s_i}$, $i\neq n$};
\end{tikzpicture}
&
\begin{tikzpicture}[scale=1.2, shift={(0,2)}]
\draw [->] (0,3) -- (1,3);
\draw [dashed] (0,2.5) -- (1,2.5);
\draw [->] (0,2) -- (1,2);
\draw [->] (0,2) -- (1,1);
\draw [->] (0,1) -- (1,2);
\draw [->] (0,1) -- (1,1);
\draw [dashed] (0,.5) -- (1,.5);
\draw [->] (0,0) -- (1,0);
\node at (.5,2.2) {++};
\node at (1,1.7) {++};
\node at (1,1.3) {++};
\node at (.5,.85) {${}--{}$};
\node at (-.3,3.05){$2n$};
\node at (1.3,3.05){$2n$};
\node at (-.5,2.05){$n+1$};
\node at (1.5,2.05){$n+1$};
\node at (-.3,1.05){$n$};
\node at (1.3,1.05){$n$};
\node at (-.3,.05){$1$};
\node at (1.3,.05){$1$};
 \node at (.5,-1){\text Diagram 2: $\psi_{s_n}$};
\end{tikzpicture}
\end{tabular}\\[1ex]	
In these two diagrams, each path from a node $k$ on the
left to a node $l$ on the right stands for an
operator on $\clh=\ell^2(\bbn)$ given as in the table. Now $\pi_{s_i}(u_{l}^{k})$ is the operator represented by the path from $k$ to $l$,
and is zero if there is no such path.
Thus, for example, $\pi_{s_i}(u_1^1)$
is $I$; $\pi_{s_i}(u_1^2)$ is zero whereas $\pi_{s_i}(u_{i+1}^{i})=-q^{N+1}$ if $i>1$.

 Next, let us explain how to represent $\pi_{s_i}\ast \pi_{s_j}$ by a diagram.
Simply keep the two diagrams representing $\pi_{s_i}$ and $\pi_{s_j}$
adjacent to each other. Identify, for each row, the node on the right side
of the diagram for $\pi_{s_i}$ with the corresponding node on the left in the diagram
for $\pi_{s_j}$. Now, $\pi_{s_i}\ast \pi_{s_j}(u_{l}^{k})$ would be
an operator on the Hilbert space $\ell^2(\bbn)\otimes \ell^2(\bbn)$ 
determined by all the paths from
the node $k$ on the left to the node $l$ on the right. It would be zero if
there is no such path and if there are more than one paths, then it would be the sum of
the operators given by each such path.
In this way, we can draw diagrams for each irreducible representation of $C(SP_q(2n))$.

 Next, we come to $\chi_{\vartheta}$. The underlying Hilbert space now is
$\ell^2(\bbz^n)\otimes \ell^2(\bbn)^{\otimes \ell(\vartheta)}$. To avoid any ambiguity, we have  explicitly
mentioned above the diagram the space on which the operator between nodes acts.
As earlier, $\chi_{\vartheta}(u_{l}^{k})$ stands for the operator
on $\ell^2(\bbz^n)\otimes \ell^2(\bbn)^{\otimes \ell(\vartheta)}$ represented by the path
from $k$ on the left to $l$ on the right.
Note that  we view $C(\bbbt^n)$ as a subalgebra of $\cll(\ell^2(\bbz^n))$.

 The following diagram  is for the representations $\chi_{\vartheta_3}$ of $C(SP_q(6))$  where $\vartheta_3 = s_3s_2s_3s_2s_1s_2s_3s_2s_1$.
 \begin{center}
 \begin{tikzpicture}[scale=1.2]
 \draw [-] (0,5) node{$\bullet$} -- (1,5) node{$\bullet$} -- (2,5) node{$\bullet$} --
     (2.5,5) node[above]{$+$} -- (3,5) node{$\bullet$} -- (4,5) node{$\bullet$}
	 -- (5,5) node{$\bullet$} -- (6,5) node{$\bullet$} -- (7,5) node{$\bullet$} 
	 -- (7.5,5) node [above]{${}+$} -- (8,5) node{$\bullet$} -- (9,5) node{$\bullet$} -- 
	 (10,5) node{$\bullet$} -- (11,5) node{$\bullet$}-- 
	 (11.5,5) node [above]{${}+$} -- (12,5) node{$\bullet$};
 \draw [-] (0,4) node{$\bullet$} -- (1,4) node{$\bullet$} -- (1.5,4) node[above]{$+$} 
       -- (2,4) node{$\bullet$} --
     (3,4) node{$\bullet$} -- (4,4) node{$\bullet$} --(4.5,4) node[above]{$+$}
	 -- (5,4) node{$\bullet$} -- (6,4) node{$\bullet$} -- (6.5,4) node[above]{$+$} -- (7,4) node{$\bullet$}  -- (7.5,4) node [below]{${}-$} -- (8,4) node{$\bullet$}  -- (8.5,4) node [above]{${}+$} -- (9,4) node{$\bullet$}  -- (10,4) node{$\bullet$}  -- (10.5,4) node [above]{${}+$} -- (11,4) node{$\bullet$}  -- (11.5,4) node [below]{${}-$} -- (12,4) node{$\bullet$} ;
 \draw [-] (0,3) node{$\bullet$} -- (.5,3) node[above]{$+$} -- (1,3) node{$\bullet$} 
       -- (2,3) node{$\bullet$} --
     (3,3) node{$\bullet$} -- (3.5,3) node[above]{$+\!+$} -- (4,3) node{$\bullet$} --
	 (4.5,3) node[below]{${}-$}
	 -- (5,3) node{$\bullet$} -- (5.5,3) node[above]{$+\!+$} -- (6,3) node{$\bullet$} -- (6.5,3) node[below]{${}-$} --(7,3) node{$\bullet$}   -- (8,3) node{$\bullet$}   -- (8.5,3) node [below]{${}-$} -- (9,3) node{$\bullet$}   -- (9.5,3) node [above]{${}+\!+$} -- (10,3) node{$\bullet$}   -- (10.5,3) node [below]{${}-$} -- (11,3) node{$\bullet$}   -- (12,3) node{$\bullet$}  ;
 \draw [-] (0,2) node{$\bullet$} -- (.5,2) node[below]{${}-$} -- (1,2) node{$\bullet$} 
       -- (2,2) node{$\bullet$} --
     (3,2) node{$\bullet$} -- (3.5,2) node[below]{$-\,-$} -- (4,2) node{$\bullet$} --
	 (4.5,2) node[above]{$+$}
	 -- (5,2) node{$\bullet$} -- (5.5,2) node[below]{$-\,-$} -- (6,2) node{$\bullet$} -- (6.5,2) node[above]{$+$} --(7,2) node{$\bullet$}   -- (8,2) node{$\bullet$}   -- (8.5,2) node [above]{${}+$} -- (9,2) node{$\bullet$}   -- (9.5,2) node [below]{${}-$} -- (10,2) node{$\bullet$}   -- (10.5,2) node [above]{${}+$} -- (11,2) node{$\bullet$}   -- (12,2) node{$\bullet$}  ;
 \draw [-] (0,1) node{$\bullet$} -- (1,1) node{$\bullet$} -- (1.5,1) node[below]{$-$} 
       -- (2,1) node{$\bullet$} --
     (3,1) node{$\bullet$} -- (4,1) node{$\bullet$} --(4.5,1) node[below]{$-$}
	 -- (5,1) node{$\bullet$} -- (6,1) node{$\bullet$} -- (6.5,1) node[below]{$-$} --(7,1) node{$\bullet$}  -- (7.5,1) node [above]{${}+$} -- (8,1) node{$\bullet$}  -- (8.5,1) node [below]{${}-$} -- (9,1) node{$\bullet$}  -- (10,1) node{$\bullet$}  -- (10.5,1) node [below]{${}-$} -- (11,1) node{$\bullet$}  -- (11.5,1) node [above]{${}+$} -- (12,1) node{$\bullet$} ;
 \draw [-] (0,0) node{$\bullet$} -- (1,0) node{$\bullet$} -- (2,0) node{$\bullet$} --
     (2.5,0) node[below]{$-$} -- (3,0) node{$\bullet$} -- (4,0) node{$\bullet$}
	 -- (5,0) node{$\bullet$} -- (6,0) node{$\bullet$} -- (7,0) node{$\bullet$}  -- (7.5,0) node [below]{${}-$} -- (8,0) node{$\bullet$}  -- (9,0) node{$\bullet$}  -- (10,0) node{$\bullet$}  -- (11,0) node{$\bullet$}  -- (11.5,0) node [below]{${}-$} -- (12,0) node{$\bullet$} ;
 \draw [-] (7,0) -- (8,1) node [below]{${}+\;$} -- (9,2) node [below]{${}+\;$} -- (10,3) node [below]{${}+\!+$} -- (11,4) node [below]{${}-\;$} -- (12,5) node [below]{${}-\;$};
 \draw [-] (7,5) -- (8,4) node [above] {${}-\;\;$} -- (9,3) node [above] {${}-\;\;$} -- (10,2) node [above] {${}+\!+$} -- (11,1) node [above] {${}+\;$} -- (12,0) node [above] {${}+\;$};
 \draw [-] (4,4) -- (5,3) node[above]{${}-$} -- (6,2) node [above]{${}+\!+$} -- (7,1)  node [above]{${}+$}-- (8,0) node [above]{${}+\;$};
 \draw [-] (4,1) -- (5,2) node[below]{${}+\;$} -- (6,3) node[below]{${}+\!+$} -- (7,4) node[below]{${}-\;$} -- (8,5) node [below]{${}-\;$};
 \draw [-] (3,2) -- (4,3) node [below]{${}+\!+$} -- (5,4) node [below]{${}-\;$};
 \draw [-] (3,3) -- (4,2) node[above]{${}+\!+$} -- (5,1) node[above]{${}+\;$};
 \draw[-] (6,1) -- (7,2) node[below]{${}+\;$};
 \draw[-] (6,4) -- (7,3) node[above]{${}-\;$};
 \draw [-] (8,2) -- (9,1) node [above]{${}+\;$};
 \draw [-] (8,3) -- (9,4) node [below]{${}-\;$};
 \draw [-] (10,1) -- (11,2) node [below]{${}+\;$};
 \draw [-] (10,4) -- (11,3) node [above]{${}-\;$};
 \draw [-] (11,0) -- (12,1) node [below]{${}+\;$};
 \draw [-] (11,5) -- (12,4) node [above]{${}-\;$};
 \node at (0.5,5.7){$\ell^2(\bbz)$};
 \node at (1.5,5.7){$\ell^2(\bbz)$};
 \node at (2.5,5.7){$\ell^2(\bbz)$};
 \node at (3.5,5.7){$\ell^2(\bbn)$};
 \node at (4.5,5.7){$\ell^2(\bbn)$};
 \node at (5.5,5.7){$\ell^2(\bbn)$};
 \node at (6.5,5.7){$\ell^2(\bbn)$};
 \node at (7.5,5.7){$\ell^2(\bbn)$};
 \node at (8.5,5.7){$\ell^2(\bbn)$};
 \node at (9.5,5.7){$\ell^2(\bbn)$};
 \node at (10.5,5.7){$\ell^2(\bbn)$};
 \node at (11.5,5.7){$\ell^2(\bbn)$};
 \node at (1,5.7){$\otimes$};
 \node at (2,5.7){$\otimes$};
 \node at (3,5.7){$\otimes$};
 \node at (4,5.7){$\otimes$};
 \node at (5,5.7){$\otimes$};
 \node at (6,5.7){$\otimes$};
 \node at (7,5.7){$\otimes$};
 \node at (8,5.7){$\otimes$};
 \node at (9,5.7){$\otimes$};
 \node at (10,5.7){$\otimes$};
 \node at (11,5.7){$\otimes$};
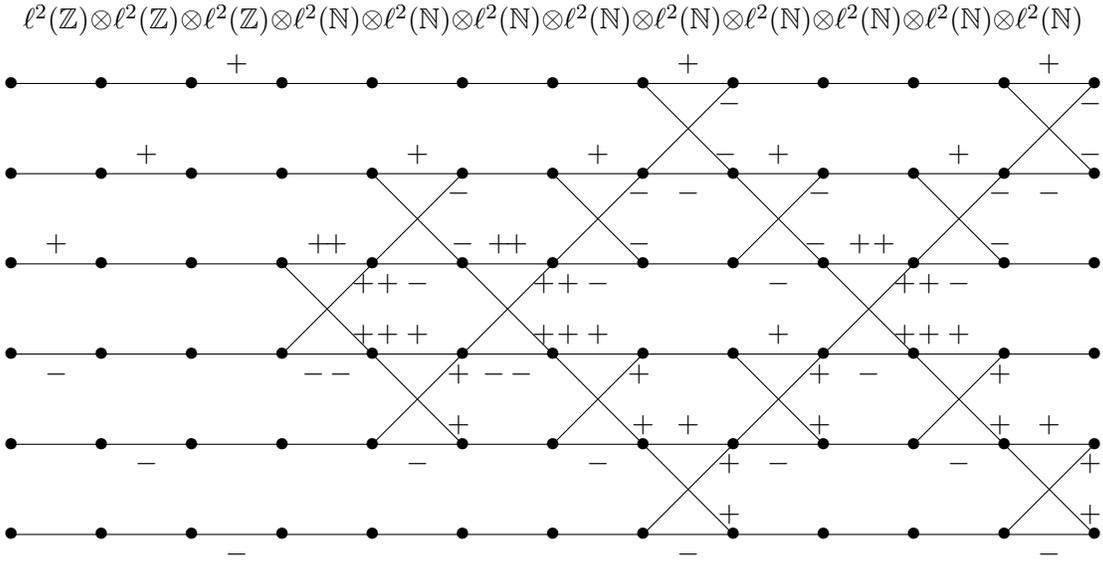
 \node at (6,-1){\text Diagram 3: $\chi_{\vartheta_3}$};
 \end{tikzpicture}
  \end{center}

 Let $\omega_n=s_1s_2\cdots s_{n-1}s_ns_{n-1}\cdots s_1$.
 The following diagram is for the representation $\pi_{\omega_3}$ of $C(SP_{q}(6))$.
 
 \begin{center}
 \begin{tikzpicture}[scale=1.2]
 \draw [-] (0,5) -- (.5,5) node [above]{${}+$} -- (1,5) -- (2,5) -- (3,5) -- (4,5)-- (4.5,5) node [above]{${}+$} -- (5,5);
 \draw [-] (0,4) -- (.5,4) node [below]{${}-$} -- (1,4) -- (1.5,4) node [above]{${}+$} -- (2,4) -- (3,4) -- (3.5,4) node [above]{${}+$} -- (4,4) -- (4.5,4) node [below]{${}-$} -- (5,4);
 \draw [-] (0,3) -- (1,3) -- (1.5,3) node [below]{${}-$} -- (2,3) -- (2.5,3) node [above]{${}+\!+$} -- (3,3) -- (3.5,3) node [below]{${}-$} -- (4,3) -- (5,3);
 \draw [-] (0,2) -- (1,2) -- (1.5,2) node [above]{${}+$} -- (2,2) -- (2.5,2) node [below]{${}-$} -- (3,2) -- (3.5,2) node [above]{${}+$} -- (4,2) -- (5,2);
 \draw [-] (0,1) -- (.5,1) node [above]{${}+$} -- (1,1) -- (1.5,1) node [below]{${}-$} -- (2,1) -- (3,1) -- (3.5,1) node [below]{${}-$} -- (4,1) -- (4.5,1) node [above]{${}+$} -- (5,1);
 \draw [-] (0,0) -- (.5,0) node [below]{${}-$} -- (1,0) -- (2,0) -- (3,0) -- (4,0) -- (4.5,0) node [below]{${}-$} -- (5,0);
 \draw [-] (0,0) -- (1,1) node [below]{${}+\;$} -- (2,2) node [below]{${}+\;$} -- (3,3) node [below]{${}+\!+$} -- (4,4) node [below]{${}-\;$} -- (5,5) node [below]{${}-\;$};
 \draw [-] (0,5) -- (1,4) node [above] {${}-\;\;$} -- (2,3) node [above] {${}-\;\;$} -- (3,2) node [above] {${}+\!+$} -- (4,1) node [above] {${}+\;$} -- (5,0) node [above] {${}+\;$};
 \draw [-] (0,1) -- (1,0) node [above]{${}+\;$};
 \draw [-] (0,4) -- (1,5) node [below]{${}-\;$};
 \draw [-] (1,2) -- (2,1) node [above]{${}+\;$};
 \draw [-] (1,3) -- (2,4) node [below]{${}-\;$};
 \draw [-] (3,1) -- (4,2) node [below]{${}+\;$};
 \draw [-] (3,4) -- (4,3) node [above]{${}-\;$};
 \draw [-] (4,0) -- (5,1) node [below]{${}+\;$};
 \draw [-] (4,5) -- (5,4) node [above]{${}-\;$};
 \node at (0,5){$\bullet$};
 \node at (1,5){$\bullet$};
 \node at (2,5){$\bullet$};
 \node at (3,5){$\bullet$};
 \node at (4,5){$\bullet$};
 \node at (5,5){$\bullet$};
 \node at (0,4){$\bullet$};
 \node at (1,4){$\bullet$};
 \node at (2,4){$\bullet$};
 \node at (3,4){$\bullet$};
 \node at (4,4){$\bullet$};
 \node at (5,4){$\bullet$};
 \node at (0,3){$\bullet$};
 \node at (1,3){$\bullet$};
 \node at (2,3){$\bullet$};
 \node at (3,3){$\bullet$};
 \node at (4,3){$\bullet$};
 \node at (5,3){$\bullet$};
 \node at (0,2){$\bullet$};
 \node at (1,2){$\bullet$};
 \node at (2,2){$\bullet$};
 \node at (3,2){$\bullet$};
 \node at (4,2){$\bullet$};
 \node at (5,2){$\bullet$};
 \node at (0,1){$\bullet$};
 \node at (1,1){$\bullet$};
 \node at (2,1){$\bullet$};
 \node at (3,1){$\bullet$};
 \node at (4,1){$\bullet$};
 \node at (5,1){$\bullet$};
 \node at (0,0){$\bullet$};
 \node at (1,0){$\bullet$};
 \node at (2,0){$\bullet$};
 \node at (3,0){$\bullet$};
 \node at (4,0){$\bullet$};
 \node at (5,0){$\bullet$};
 \node at (.5,5.7){$\ell^2(\bbn)$};
 \node at (1.5,5.7){$\ell^2(\bbn)$};
 \node at (2.5,5.7){$\ell^2(\bbn)$};
 \node at (3.5,5.7){$\ell^2(\bbn)$};
 \node at (4.5,5.7){$\ell^2(\bbn)$};
 \node at (1,5.7){$\otimes$};
 \node at (2,5.7){$\otimes$};
 \node at (3,5.7){$\otimes$};
 \node at (4,5.7){$\otimes$};
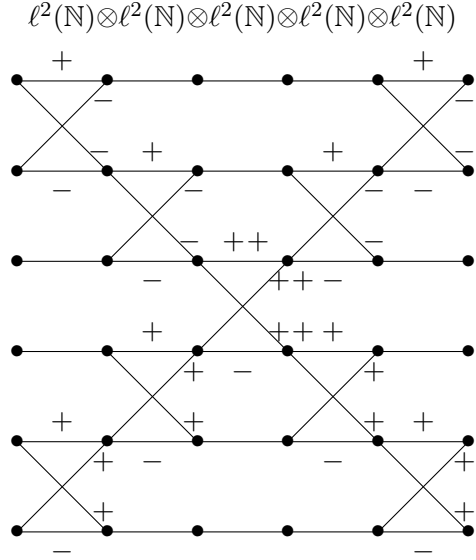
 \node at (2.5,-1){\text Diagram 4: $\pi_{\omega_3}$};
 \end{tikzpicture}
 \end{center}

\newsection{The quotient space $C(SP_{q}(2n)/SP_{q}(2n-2))$}
In this section, we  recall quantum homogeneous space $C(SP_{q}(2n)/SP_{q}(2n-2))$ and show that it is the 
$C^*$- subalgebra of $C(SP_q(2n))$ generated by $\left\{u_m^1,u_m^{2n} : m \in \{1,2,\cdots 2n\}\right\}$.

The Weyl group $W_{n-1}$ of $sp_{2n-2}$ can  be realized as a subgroup of the Weyl group $W_{n}$ of $sp_{2n}$ 
generated by ${s_{2},s_{3},\cdots s_{n}}$ and hence, the longest word $\vartheta_{n-1}$ in $W_{n-1}$ is a
subword of the longest word $\vartheta_{n}$ in $W_{n}$ which can easily be seen from proposition $\ref{pw1}$. 
This shows that $C(SP_{q}(2n-2))$ is a subgroup of $C(SP_{q}(2n))$, i.e. there is a $C^{*}$-epimorphism 
$\phi:C(SP_{q}(2n))\rightarrow C(SP_{q}(2n-2))$ obeying $\Delta\phi=(\phi\otimes \phi)\Delta$. More precisely, 
let $\sigma:\scrt \rightarrow \bbc$ is the homomorphism for which $\sigma(S) = 1$. Define $\phi$ to be the restriction 
of $1^{\otimes n-1}\otimes ev_{1}\otimes 1^{\otimes (n-1)^{2}}\otimes \sigma^{\otimes (2n-1)}$ to $\chi_{\vartheta_{n}}(C(SP_{q}(2n)))$
which is contained in $C(T^{n})\otimes \scrt^{\otimes n^{2}}$. Here $ev_{1}$ denote the evaluation map at
1 i.e. $ev_{1}:C(T) \rightarrow \bbc$ such that $ev_{1}(f)=f(1)$.  Image of $\phi$ is equal to $\chi_{\vartheta_{n-1}}(C(SP_{q}(2n-2)))$ as,
\[
\phi(\chi_{w_{n}}(u_{j}^{i}))=\begin{cases}
                             \chi_{\vartheta_{n-1}}(v_{j}^{i}), 
                             & \mbox{ if } i \neq 1 \mbox{ or } 2n, \mbox{ or } j \neq 1 \mbox{ or } 2n, \cr
			 \delta_{ij}, & \mbox{ otherwise. } \cr
			\end{cases}
\] 
where $v_{j}^{i}$ are generators of $C(SP_{q}(2n-2))$. In such a case, one defines the quotient space $C(SP_{q}(2n)/SP_{q}(2n-2))$ by,
\[
         C(SP_{q}(2n)/SP_{q}(2n-2)) = \left\{a\in C(SP_{q}(2n)) : (\phi\otimes id)\Delta(a) = I\otimes a\right\}. 
\] 
\bthm \label{quotient}
The quotient space $C(SP_{q}(2n)/SP_{q}(2n-2))$ is the $C^*$-algebra generated by $\left\{u_m^1,u_m^{2n} : m \in \{1,2,\cdots 2n\}\right\}$.
\ethm
We first prove one proposition which will be needed in the proof of above theorem. For that, define 
$\clb$ to be the $*$-algebra generated by $\left\{u_m^1,u_m^{2n} : m \in \{1,2,\cdots 2n\}\right\}$. 
Consider $\clb$ as a $U_q(sp(2n))$-module with the following action:
\[
 f(a)=(1\otimes \langle f,.\rangle)\Delta a.
\]
where $f \in U_q(sp(2n)), a \in \clb$ and $\langle.,.\rangle$ is the pairing given in theorem \ref{pairing}. Fix a $n$-tuple of integers of the form
$(r,s,0,\cdots,0)$ satisfying $r \geq s$. 
We call $b \in \clb$ highest weight vector  with highest weight $(r,s,0.\cdots,0)$ if 
\begin{IEEEeqnarray}{rCll}
K_1(b)&=&q^{r-s}b,\nonumber \\
K_2(b)&=&q^{s}b, \nonumber \\
K_i(b)&=& b  &\mbox{for all} \quad i\geq 2, \nonumber \\
E_i(b)&=&0  &\mbox{for all} \quad i \in \left\{1,\cdots , 2n\right\}. \nonumber 
\end{IEEEeqnarray}
\bppsn \label{hwv}
There exist $r-s+1$ linearly independent highest weight vectors in $\clb$ with highest weight $(r,s,0,\cdots,0)$ for all $r,s \in \bbn$ satisfying 
$r\geq s$.
\eppsn
\prf
Let $x= u_{2n-1}^1,y=u_{2n-1}^{2n},z=u_{2n}^1$ and $w=u_{2n}^{2n}$. Let $p=r-s$.  It is easy to see that
\begin{displaymath}
 E_i(x)=E_i(y)=E_i(z)=E_i(w)=0 \quad \mbox{for}\quad  i >1.
\end{displaymath}
 Also, 
 \begin{displaymath}
  E_1(x)=-z, E_1(y)=-w, E_1(z)=E_2(w)=0.
 \end{displaymath}
 Further, 
 \begin{displaymath}
  K_1(x)=q^{-1}x, K_1(y)=q^{-1}y, K_1(z)=qz, K_1(w)=qw, 
 \end{displaymath}
 \begin{displaymath}
  K_2(x)=qx,K_2(y)=qy, K_2(z)=z, K_2(w)=w.
 \end{displaymath}

 and for $i>2$, $K_i$ fixes these elements. Now, using the relations $(\ref{relations})$ and  
 the facts that $\Delta(E_1)=E_1\otimes K_1+1\otimes E_1$ 
and $\Delta(K_i)=K_i\otimes K_i$, we get
\begin{IEEEeqnarray}{rCl}
 E_1(y^sz^s)&=&C_1^0y^{s-1}z^sw. \nonumber \\
 E_1(xy^{s-1}z^{s-1}w)&=&C_1^1y^{s-1}z^sw+C_2^1xy^{s-2}z^{s-1}w. \nonumber \\
 \cdots & & \cdots \nonumber\\
 E_1(x^{s-1}yzw^{s-1})&=&C_1^{s-1}x^{s-2}yz^2w^{s-1}+C_2^{s-1}x^{s-1}zw^{s}.\nonumber \\
 E_1(x^sw^s)&=&C_1^sx^{s-1}zw^s \nonumber
\end{IEEEeqnarray}
where $C_i^j$'s are nonzero constants. This shows that we can choose nonzero constants $c_1 \cdots c_s$ such that 
$E_1(y^sz^s+c_1xy^{s-1}z^{s-1}w+\cdots +c_sx^sw^s)=0$. Let $\omega=s_1s_2\cdots s_{n-1}s_ns_{n-1}\cdots s_1$ and 
$\omega^{'}=s_1s_2\cdots s_{n-1}s_ns_{n-1}\cdots s_2$. Let $\pi_{\omega}$ and $\pi_{\omega^{'}}$
be the representations of $C(SP_q(2n))$ as defined in subsection $(2.2)$. It is easy to see that $\pi_{\omega^{'}}(z)=0$. Hence,
\[
\pi_{\omega^{'}}(y^sz^s+c_1xy^{s-1}z^{s-1}w+\cdots +c_sx^sw^s)(e_0\otimes e_0\otimes \cdots \otimes e_0)=
\pi_{\omega^{'}}(c_sx^sw^s)(e_0\otimes e_0\otimes \cdots \otimes e_0)\neq 0
\]
which shows that $y^sz^s+c_1xy^{s-1}z^{s-1}w+\cdots +c_sx^sw^s \neq 0$. Since $\omega^{'}$ is a subword of $\omega$, $\pi_{\omega^{'}}$
factors through $\pi_{\omega}$ which implies that $\pi_{\omega}(y^sz^s+c_1xy^{s-1}z^{s-1}w+\cdots +c_sx^sw^s) \neq 0$.

Define $b_j=z^jw^{p-j}(y^sz^s+c_1xy^{s-1}z^{s-1}w+\cdots +c_sx^sw^s)$ for $j \in \left\{0,\cdots,p\right\}$. Clearly, $b_j \in \clb$.
One can directly verify that $b_j$ are elements with highest weight $(r,s,0.\cdots,0)$.  
Now, look at the $(n)^{\mbox{th}}$ position of  $\pi_{\omega}(b_j)(e_0\otimes e_0\otimes \cdots e_0)$.
One term has $e_{p-j}$ at $(n)^{\mbox{th}}$ position and other terms have $e_{\ell}$ at $(n)^{\mbox{th}}$ position where $\ell < p-j$
(see diagram 4). This proves
that $b_j$ are linearly independent.\qed 

\prf (of theorem \ref{quotient})
One can easily check that  $u_{m}^{2n}$ and  $u_{m}^{1}=\epsilon_{m}q^{\rho_{1}-\rho_{m}}(u_{2n-m+1}^{2n})^{*}$ are in $C(SP_{q}(2n)/SP_{q}(2n-2))$
for $m=1,2,\cdots 2n$. So, 
\[
C(SP_{q}(2n)/SP_{q}(2n-2)) \supseteq C^*\left\{u_m^1,u_m^{2n} : j \in \{1,2,\cdots 2n\}\right\}.
\]
 To show the equality,
consider the co-multiplication action on $C(SP_{q}(2n)/SP_{q}(2n-2))$ by the compact quantum group $C(SP_q(2n))$, 
\begin{IEEEeqnarray}{rCl}
 C(SP_{q}(2n)/SP_{q}(2n-2) &\longrightarrow & C(SP_{q}(2n)/SP_{q}(2n-2) \otimes C(SP_q(2n)) \nonumber \\
 a &\longmapsto & \Delta a \nonumber 
\end{IEEEeqnarray}
By theorem $1.5$, Podles \cite{Pod-1995aa}, we get,
\begin{IEEEeqnarray}{rCl}
C(SP_{q}(2n)/SP_{q}(2n-2)= \overline{\oplus_{\lambda \in \widehat{SP(2n)}}\oplus_{i \in I_{\lambda} }W_{\lambda,i}} \nonumber 
\end{IEEEeqnarray}
where $\lambda$ represents a finite-dimensional irreducible co-representation $u^{\lambda}$ of $C(SP_q(2n))$, $ W_{\lambda, i}$ corresponds to 
$u^{\lambda}$ for all $i \in I_{\lambda}$ and $I_{\lambda}$ is the multiplicity of $u^{\lambda}$.
Define 
\[
\cla= \oplus_{\lambda \in \widehat{SP(2n)}}\oplus_{i \in I_{\lambda} }W_{\lambda,i}.
\]
We will prove that $\cla \subseteq 
C^*\left\{u_m^1,u_m^{2n} : m \in \{1,2,\cdots 2n\}\right\}$ which will suffice to show the claim.

The finite-dimensional irreducible co-representations of $C(SP_q(2n))$ or equivalently irreducible representations of $U_q(sp(2n))$ are in a 
one-to-one correspondence with $n$-tuples of integers $\lambda=(\lambda_1,\cdots,\lambda_n)$ satisfying the inequalities 
\[
 \lambda_1 \geq \lambda_2 \geq \cdots \geq \lambda_n \geq 0
\]
Such an $n$-tuple $\lambda$ is called the highest weight of the corresponding representation which we denote by $V(\lambda)$. The restriction
of $V(\lambda)$ to the subalgebra $U_q(sp(2n-2))$ is isomorphic to a direct sum of irreducible finite-dimensional representations $V^{'}(\mu)$,
$\mu =(\mu_1,\cdots , \mu_{n-1})$ of $U_q(sp(2n-2))$ with certain multiplicity $n_{\lambda}(\mu)$. The multiplicity $n_{\lambda}(\mu)$ is equal
to the number of $n$-tuples of integers $(\nu_1, \cdots \nu_n)$ satisfying the inequalities,
\begin{IEEEeqnarray}{rCl}
 \lambda_1 \geq \nu_1 \geq \lambda_2 \geq \nu_2 \geq \cdots \geq \lambda_n \geq \nu_n \geq 0. \nonumber \\
 \nu_1 \geq \mu_1 \geq \nu_2 \geq \mu_2 \cdots  \geq \mu_{n-1} \geq \nu_n \geq 0 \nonumber. 
\end{IEEEeqnarray}
We refer to Zhelobenko~\cite{Zhe-1962aa} for more detail. 
Now, this shows that a finite-dimensional irreducible representation of $U_q(sp(2n))$ with highest weight $\lambda$ when restricted to the sublagebra 
$U_q(sp(2n-2))$ contains trivial representation if and only if $\lambda_i=0$ for all $i \geq 3$ and for such $\lambda$, the multiplicity of 
trivial representation denoted by $n_{\lambda}(0)=\lambda_1-\lambda_2+1$. By theorem $1.7$, Podles \cite{Pod-1995aa},
\begin{displaymath}
I_{\lambda}=n_{\lambda}(0)= \begin{cases}
                             \lambda_1-\lambda_2+1, & \quad \mbox{if} \quad \lambda_i=0 \quad \mbox{for all}\quad  i \geq 3, \cr
                             0, & \quad \mbox{otherwise}, \cr
                            \end{cases}
\end{displaymath}
It follows from proposition \ref{hwv} that $\clb \subseteq C^*\left\{u_m^1,u_m^{2n} : m \in \{1,2,\cdots 2n\}\right\}$ contains 
$\lambda_1-\lambda_2+1$ linearly independent highest weight vector with highest weight $(\lambda_1,\lambda_2,0,\cdots,0)$.This proves that for 
each co-representation $\lambda$ of $SP(2n)$,
$\oplus_{i\in I_{\lambda}}W_{\lambda,i} \subseteq C^*\left\{u_m^1,u_m^{2n} : m \in \{1,2,\cdots 2n\}\right\}$ which further shows that 
$\cla \subseteq C^*\left\{u_m^1,u_m^{2n} : m \in \{1,2,\cdots 2n\}\right\}$.
This proves the claim. 
\qed

 \newsection{ Quantum quaternion sphere}
 Our main aim in this section is to describe $C(SP_q(2n)/SP_q(2n-2))$  as a universal $C^*$-algebra given by a 
 finite set of generators 
 and relations.
 
\bdfn 
 We define $C^{*}$-algebra $C(H_{q}^{2n})$ of continuous functions on the quantum quaternion sphere  as the universal  
 $C^{*}$-algebra generated by elements $z_{1}$, $z_{2}$, ....$z_{2n}$ satisfying the following relations
\begin{align}
z_{i}z_{j} &=  qz_{j}z_{i} & \text{for } & i > j,  i+j \neq 2n+1 \label{c1}\\
z_{i}z_{i^{'}}  &= q^2z_{i^{'}}z_{i} -(1-q^{2})\sum_{k > i} q^{i-k}z_{k}z_{k^{'}} & \text{for } &  i > n \label{c2}\\
z_{i}^{*}z_{i^{'}} &= q^{2}z_{i^{'}}z_{i}^{*} \label{c3}\\
z_{i}^{*}z_{j} &= qz_{j}z_{i}^{*} & \text{for }  & i+j>2n+1,  i \neq j \label{c4}\\
z_{i}^{*}z_{j} &=qz_{j}z_{i}^{*} + (1-q^{2})\epsilon_{i}\epsilon_{j}q^{\rho_{i}+\rho_{j}}z_{i^{'}}z_{j^{'}}^{*} 
& \text{for }  & i+j<2n+1,  i \neq j \label{c5}\\
z_{i}^{*}z_{i} &=z_{i}z_{i}^{*} + (1-q^{2})\sum_{k>i}z_{k}z_{k}^{*} & \text{for } & i >n \label{c6}\\
z_{i}^{*}z_{i} &=z_{i}z_{i}^{*} +(1-q^{2})q^{2\rho_{i}}z_{i^{'}}z_{i^{'}}^{*}+ (1-q^{2})\sum_{k>i}z_{k}z_{k}^{*} 
   & \text{for }  & i \leq n \label{c7}\\
\sum_{i=1}^{2n}z_{i}z_{i}^{*} &= 1 \label{c8}
\end{align} 
\edfn
In what follows, we will find a faithful realization of this $C^*$-algebra on a Hilbert space.  For this, we will first find 
all irreducible representations of the above $C^*$-algebra.

 It follows from the commutation relations that $ \left\|z_{i}\right\| \leq 1, \mbox{ for } 1 \leq i \leq 2n$ and $z_{2n}$ is normal. 
 We denote $z_{2n}^{*}z_{2n}$ by $\omega$. Using the  relations $(\ref{c1}), (\ref{c3})$ and $(\ref{c5})$, we have
\begin{IEEEeqnarray}{lCl}
z_{i}\omega &=& q^{-2}\omega z_{i}, \qquad z_{i}^{*}\omega = q^{2}\omega z_{i}^{*} \quad \mbox{ for all } i \not\in \{ 1 , 2n\}, 
\label{chap5-eqn-1} \\
z_{1}\omega &=& q^{-4}\omega z_{1}, \qquad z_{1}^{*}\omega = q^{4}\omega z_{1}^{*}. \label{chap5-eqn-2}
\end{IEEEeqnarray}

\bppsn \label{chap5-ppsn-1}
Let $\pi$ be a representation of $C(H_{q}^{2n})$. Then one has 
\begin{enumerate}
\item 
$\pi(\omega) = I \qquad \mbox{ on } \;\bigcap_{i=1}^{2n-1} \ker \pi(z_{i}^{*})$,
\item 
$1_{(q^{2m+2}, q^{2m})}(\pi(\omega)) = 0 \quad \forall m \in \bbn$,
\item
$\ker(\pi(z_{i}))\subseteq \ker(\pi(z_{k}^{*})) \mbox{ for } k\geq i \mbox{ and } 1 \leq i \leq 2n$,
\item
if $u$ is a nonzero eigenvector of $\pi(\omega)$ corresponding to the eigenvalue $q^{2m}$, then 
$u \notin \ker\pi(z_{i})$  for  $1 \leq i \leq 2n-1$,
\item
either $\sigma (\pi(\omega)) = \left\{q^{2m}: m \in \bbn \right\} \bigcup \left\{0\right\}$ or $\sigma (\pi(\omega)) = \left\{0\right\}$. 
\end{enumerate}
\eppsn 
\prf
\begin{enumerate}
\item
Easy to see from $(\ref{c8})$.
\item
From the commutation relations, it follows that $z_1^*f(\omega)=f(q^4\omega)z_1^*$ and 
$z_i^*f(\omega)=f(q^2\omega)z_i^*$ for all $i \neq 1$ for all
continuous functions $f$ and hence for all $L_\infty$ functions.
Thus 
\begin{align*}
         \pi(z_1)^*1_{(q^{2n+2},q^{2n})}(\pi(\omega)) 
		 & =1_{(q^{2n+2},q^{2n})}(q^4\pi(\omega))\pi(z_1)^*\\
		 &=1_{(q^{2n-2},q^{2n-4})}(\pi(\omega))\pi(z_1)^*,\\
         \pi(z_i)^*1_{(q^{2n+2},q^{2n})}(\pi(\omega))
                      &=1_{(q^{2n+2},q^{2n})}(q^2\pi(\omega))\pi(z_i)^*\\
                       &=1_{(q^{2n},q^{2n-2})}(\pi(\omega))\pi(z_1)^*.
\end{align*}
By repeated application and using $(\ref{c8})$ and the fact that
$\sigma(\omega)\subseteq[0,1]$, it follows that $1_{(q^{2n+2},q^{2n})}(\pi(\omega))=0$.
\item
Let $h \in $ ker$(\pi(z_{i}))$ and $ i > n$. Using $(\ref{c6})$, we have
\[
\left\langle z_i^*z_ih, h \right\rangle = \left\langle z_iz_i^*h + (1-q^2)\sum_{k>i}z_kz_k^*h, h \right\rangle.
\]
Therefore it follows that
\[
\left\|z_i^*h\right\|^2 + (1-q^2)\sum_{k>i}\left\|z_k^*h\right\|= 0.
\]
Hence $\left\|z_k^*h\right\| = 0$ for all  $k\geq i$, which means 
 $h  \in  \mbox{ ker }\pi(z_{k}^{*})$ or all  $k\geq i$. 
For $i \leq n$, use $(\ref{c7})$ and follow similar steps. 

\item
From part $3$, we have ker$(z_{i})\subseteq \ker(z_{2n}^*) = \ker(z_{2n})= \ker(\omega)$.
Now if $u$ is a non-zero eigenvector of $\pi(\omega)$ corresponding to eigenvalue $q^{2m}$ for
some $m \in \bbn$, then $u \notin \ker(z_{2n}^*)$. Hence $u \notin \ker(z_{i})$ for $1 \leq i \leq 2n$. 
\item
From part $2$ and the fact that $\left\|\omega\right\|\leq 1$, it follows that 
$\sigma(\pi(\omega))\subseteq\left\{q^{2m}: m \in \bbn\right\} \bigcup \left\{0\right\}$. Define
\[
 A = \left\{m \in \bbn : q^{2m} \in  \sigma(\pi(\omega))\right\}.
\]
If $A = \emptyset$, we have $\sigma(\pi(\omega)) = \left\{0\right\}$. If $A \neq \emptyset$, define
\[ 
m_{0}=\mbox{ inf }\left\{m \in \bbn : q^{2m} \in  \sigma(\pi(\omega))\right\}.
\] 
Let $u$ be a nonzero eigenvector corresponding to $q^{2m_{0}}$. Assume $u \notin  \ker \pi(z_{i}^{*})$
for some $i \in \left\{1,2,\cdots ,2n-1\right\}$. Then from $(\ref{c1})$, it follows that  $\pi(z_{i}^{*})u$ 
is a nonzero eigenvector corresponding to the eigenvalue $q^{2m_{0}-2}$ or $q^{2m_{0}-4}$ depending on whether 
$i \neq 1$ or $i=1$, which contradicts the fact that $m_{0}$ is $\inf A$. Hence $u \in \bigcap_{i=1}^{2n-1}\ker \pi(z_{i}^{*})$.
As  $\pi(\omega) = I$ on $\bigcap_{i=1}^{2n-1}\ker \pi(z_{i}^{*})$, we get $m_{0} = 0$. From part $4$, it follows that 
$u \notin \ker \pi(z_{i})$ for any $i \in \left\{1,2,\cdots 2n\right\}$. Now applying $(\ref{c1})$, we have
$\pi(z_{2}^{m})u$ is a nonzero eigenvector corresponding to eigenvalue $q^{2m} \mbox{ for all } m \in \bbn$. This proves the claim. 
\end{enumerate} 
\qed

 Let $\pi$ be a representation of $C(H_{q}^{2n})$ in a Hilbert space $\clh$. From $(\ref{c1})$, it follows that $\ker(\pi(\omega))$ is 
 an invariant subspace for $\pi$. Therefore if $\pi$ is irreducible, then either $\pi(\omega)=0$ or $\ker(\pi(\omega))=0$.
 Assume $\pi(\omega) \neq 0$. Then ker$(\pi(\omega))=0$, and by part $5$, we have
 $\sigma (\pi(\omega)) = \left\{q^{2m}: m \in \bbn \right\} \bigcup \left\{0\right\}$. Hence $\clh$ decomposes as
\[\clh = \oplus_{m \in \bbn}\clh_{m}. 
\] 
where $\clh_{m}$ is the eigenspace of $\pi(\omega)$ corresponding to the eigenvalue $q^{2m}$.
It is clear from $(\ref{c1})$ that for $1 < i < 2n$, the operator $\pi(z_{i})$ sends $\clh_{m}$ into 
$\clh_{m+1}$ and $\pi(z_{i}^{*})$ sends $\clh_{m}$ into $\clh_{m-1}$, $\pi(z_{1})$ sends 
$\clh_{m}$ into $\clh_{m+2}$ and $\pi(z_{1}^{*})$ sends $\clh_{m}$ into $\clh_{m-2}$.
Also, both $\pi(z_{2n})$ and $\pi(z_{2n}^{*})$ keep $\clh_{m}$ invariant. Observe that $\pi(z_{2n})|_{\clh_{0}}$ is an unitary operator. 
\bppsn
Let $u \in \bigcap_{i=1}^{2n-1}$ ker $\pi(z_{i}^{*})$. Then
\[
\pi(z_{2n})u \in \bigcap_{i=1}^{2n-1}  \ker \pi(z_{i}^{*}),\qquad \qquad \pi(z_{2n}^{*})u \in \bigcap_{i=1}^{2n-1} \ker \pi(z_{i}^{*}). 
\]
\eppsn 
\prf
 We need to show that 
       $\pi(z_{i})\pi(z_{2n})u= \pi(z_{i}^{*})\pi(z_{2n})u=0$  for all  
	   $i \in \left\{1,2,\cdots 2n-1\right\}$, 
 which  follows from $(\ref{c1}),(\ref{c3})$ and $(\ref{c4})$. 
 \qed

Let $K$ be a subspace of $\cap_{i=1}^{2n-1} \ker \pi(z_i^*)$ invariant under the $C^*$-algebra generated by $\pi(z_{2n})$. Define
\begin{displaymath}
\clh^{K} = \mbox{ linear span } \bigg\{ \pi(z_{1})^{\alpha_{1}}\pi(z_{2})^{\alpha_{2}}\cdots \pi(z_{2n-1})^{\alpha_{2n-1}}h: h \in K \bigg\}.
\end{displaymath}

\blmma \label{invariant}
Let $\pi$ be an irreducible representation of $C(H_q^{2n})$ such that $\pi(z_{2n}) \neq 0$. Then
$\clh^{K}$ is an invariant subspace of $\pi$.
\elmma 
\prf
Let  $h \in K$. Define 
\[ 
h(\alpha_{1},\alpha_{2},\cdots ,\alpha_{2n-1})
= \pi(z_{1})^{\alpha_{1}}\pi(z_{2})^{\alpha_{2}}\cdots \pi(z_{2n-1})^{\alpha_{2n-1}}h.
\]
It is clear that $\pi(z_{2n})$ keeps $\clh^{K}$ invariant, as
\[\pi(z_{2n})h(\alpha_{1},\alpha_{2},\cdots ,\alpha_{2n-1}) 
= q^{(\sum_{l=1}^{2n-1}\alpha_{l})+\alpha_{1}}h(\alpha_{1},\alpha_{2},\cdots ,\alpha_{2n-1}).
\]
 For $1 \leq i \leq n$,
\[
\pi(z_{i})h(\alpha_{1},\alpha_{2},\cdots ,\alpha_{2n-1})
= q^{\sum_{l=1}^{i-1}\alpha_{l}}h(\alpha_{1},\cdots ,\alpha_{i-1},\alpha_{i}+1,\alpha_{i+1},\cdots ,\alpha_{2n-1}) \in \clh^{K}.
\]
For $i=2n-1$,
\[
\pi(z_{2n-1})h(\alpha_{1},\alpha_{2},\cdots ,\alpha_{2n-1})
= q^{\alpha_{1}}\pi(z_{1})^{\alpha_{1}}\pi(z_{2n-1})\pi(z_{2})^{\alpha_{2}}\pi(z_{3})^{\alpha_{3}}\cdots \pi(z_{2n-1})^{\alpha_{2n-1}}h.
\] 
Repeated application of $(\ref{c2})$ gives
\[
z_iz_{i^{'}}^{m} = q^{2m}z_{i^{'}}^{m}z_i-(1-q^{2m})\sum_{k>i}q^{i-k}z_{i^{'}}^{m-1}z_kz_{k^{'}}.
\]
Hence we have 
\begin{IEEEeqnarray}{lCl}
         \IEEEeqnarraymulticol{3}{l}{\pi(z_{2n-1})
		 h(\alpha_{1},\alpha_{2},\cdots ,\alpha_{2n-1})}\nonumber \\
         &=& q^{\alpha_{1}}\pi(z_{1})^{\alpha_{1}}
		 \pi(q^{2\alpha_2}z_2^{\alpha_2}z_{n-1}-q(1-q^{2\alpha_2})
                     z_2^{\alpha_2-1}z_{2n}z_1)\pi(z_{3})^{\alpha_{3}}\cdots 
		     \pi(z_{2n-1})^{\alpha_{2n-1}}h.\nonumber\\
          &=& q^{\alpha_{1}+2\alpha_2+\sum_{l=3}^{2n-3}}
		  h(\alpha_{1},\alpha_{2},\cdots ,\alpha_{2n-1}+1) \nonumber\\
          &&\>   -q^{\alpha_1+1}(1-q^{2\alpha_2}) \pi(z_{1})^{\alpha_{1}}
		  \pi(z_2^{\alpha_2-1}z_{2n}z_1)h(0,0,\alpha_{3},\cdots ,\alpha_{2n-1}).\nonumber
\end{IEEEeqnarray} 
We have shown above that $\pi(z_{1}),\pi(z_{2})$ and $\pi(z_{2n})$ keep $\clh^{K}$ invariant. Hence
\[
\pi(z_{2n-1})h(\alpha_{1},\alpha_{2},\cdots ,\alpha_{2n-1}) \in \clh^K.
\]
Similarly, by using backward induction  we can show that $\clh^{K}$ is invariant under
the actions of $\pi(z_1),\pi(z_2),\cdots \pi(z_n)$. Also, we have 
\[
\pi(z_{2n}^*)h(\alpha_{1},\alpha_{2},\cdots ,\alpha_{2n-1})=
q^{(\sum_{l=1}^{2n-1}\alpha_{l})+\alpha_{1}}h(\alpha_{1},\alpha_{2},\cdots ,\alpha_{2n-1}).
\]
This shows that $\pi(z_{2n}^*)$ keeps $\clh^{K}$ invariant. By applying $(\ref{c5})$ and $(\ref{c6})$  repeatedly, we get
\begin{align*}
       z_{2n-1}^{*}z_{1}^m  &=   q^mz_{1}^{m}z_{2n-1}^{*} +  
	mq^{m}(1-q^{2})\epsilon_{2n-1}\epsilon_{1}q^{\rho_{2n-1}+
	\rho_{1}}z_1^{m-1}z_{2}z_{2n}^{*},\\
       z_{2n-1}^*z_{2n-1}^m  &=  z_{2n-1}^mz_{2n-1}^*+
	   (1-q^{2m})z_{2n-1}^{m-1}\omega.
\end{align*}
Hence we have 
\begin{IEEEeqnarray}{lCl}
    \IEEEeqnarraymulticol{3}{l}{\pi(z_{2n-1}^*)
	h(\alpha_{1},\alpha_{2},\cdots ,\alpha_{2n-1})}\nonumber \\
     &=& q^{\alpha_1}\pi(z_{1}^{\alpha_1}z_{2n-1}^{*} + 
	 \alpha_1 q^{\alpha_1}(1-q^{2})\epsilon_{2n-1}
           \epsilon_{1}q^{\rho_{2n-1}+\rho_{1}}z_1^{\alpha_1 -1}
		   z_{2}z_{2n}^{*})\pi(z_{2})^{\alpha_{2}}\cdots  
		      \pi(z_{2n-1})^{\alpha_{2n-1}}h\nonumber\\
      &=& q^{(\sum_{l=1}^{2n-2} l) + \alpha_2}
	  \pi(z_{1})^{\alpha_{1}}\cdots \pi(z_{2n-2})^{\alpha_{2n-2}}
            \pi(z_{2n-1})^*\pi(z_{2n-1})^{\alpha_{2n-1}}h\nonumber \\
       &&\> -\alpha_1 q^{\alpha_1}(1-q^{2})q^{\rho_{2n-1}+\rho_{1}}
	      \pi(z_1^{\alpha_1 -1})\pi(z_{2})\pi(z_{2n}^{*})\pi(z_{2})^{\alpha_{2}}
		  \cdots \pi(z_{2n-1})^{\alpha_{2n-1}}h\nonumber\\
       &=& (1-q^{2\alpha _1})q^{(\sum_{l=1}^{2n-2} l) +
	    \alpha_2}\pi(z_{1})^{\alpha_{1}}\cdots \pi(z_{2n-2})^{\alpha_{2n-2}}
		\pi(z_{2n-1})^{\alpha_{2n-1}-1}\pi(\omega)h\nonumber \\
       &&\> -\alpha_1 q^{\alpha_1}(1-q^{2})q^{\rho_{2n-1}+\rho_{1}}
	   \pi(z_1^{\alpha_1 -1})\pi(z_{2})\pi(z_{2n}^{*})\pi(z_{2})^{\alpha_{2}}
	   \cdots  \pi(z_{2n-1})^{\alpha_{2n-1}}h.\nonumber
\end{IEEEeqnarray}
Since $\pi(z_{1}),\pi(z_{2})$ and $\pi(z_{2n}^*)$ keep $\clh^{K}$ invariant, $\clh^{K}$ is invariant under the action of $\pi(z_{2n-1}^*)$.
By using backward induction and following similar steps, we can show that $\clh^{K}$ is invariant for $\pi$. 
\qed  

It follows from the lemma above that if $K$ is an invariant subspace for $\bigcap_{i=1}^{2n-1}  \ker \pi(z_{i}^{*})$, then
$\clh^{K}$ is an invariant subspace for $\pi$ and is a proper invariant subspace for $\pi$ if $K$ is a proper subspace of
$\bigcap_{i=1}^{2n-1}  \ker \pi(z_{i}^{*})$. Therefore, if $\pi$ is an irreducible representation, then the space 
$\bigcap_{i=1}^{2n-1} \ker \pi(z_{i}^{*})$ is one dimensional.
\bcrlre \label{cr2}
Let $\pi$ be an irreducible representation such that $\pi(z_{2n}) \neq 0$. Let $u$ be a unit vector in $\bigcap_{i=1}^{2n-1} \ker \pi(z_{i}^{*})$.
Then
\begin{displaymath}
\clh_{m} = \mbox{ linear span }\bigg\{\ \pi(z_{1})^{\alpha_{1}}\pi(z_{2})^{\alpha_{2}}\cdots \pi(z_{2n-1})^{\alpha_{2n-1}}u : 
\Big(\sum_{l=1}^{2n-1}\alpha_{l}\Big)+\alpha_{1}=m \bigg\}.
\end{displaymath}
\ecrlre
\prf 
It follows from Lemma \ref{invariant},  equation $(\ref{chap5-eqn-1})$ and equation $(\ref{chap5-eqn-2})$.
\qed

Define
\begin{align} \label{chap5-eqn-u}
      u_{\alpha_{2},\cdots ,\alpha_{2n-1},\alpha_0} 
      &= \pi(z_{2n-1})^{\alpha_{2n-1}}\pi(z_{2n-2})^{\alpha_{2n-2}}
             \cdots \pi(z_{2})^{\alpha_{2}}[\pi(z_{n}),\pi(z_{n+1})]^{\alpha_0}u,
\end{align} 
where $\alpha_{i} \in \bbn$ and $u$ is as in Corollary \ref{cr2}. 
Now we develop some tools by analyzing the defining relations more closely.  
\bppsn \label{chap5-ppsn-2}
Let $\pi$ be an irreducible representation of $C(H_{q}^{2n})$ such that 
$\pi(z_{2n}) \neq 0$. Then
\begin{enumerate}
\item
for $i > n$,
\begin{displaymath}
\pi(z_{i})^{*}\pi(z_{i})^{m} = \pi(z_{i})^{m}\pi(z_{i})^{*} + (1-q^{2m})\sum_{k>i}\pi(z_{i})^{m-1}\pi(z_{k})\pi(z_{k})^{*}; 
\end{displaymath}
\item
for $i \leq n$,
\begin{IEEEeqnarray*}{rCl}
      \pi(z_{i})^{*}\pi(z_{i})^{m}   
	  &=& \pi(z_{i})^{m}\pi(z_{i})^{*} +q^{2\rho_{i}}(1-q^{2m})\pi(z_{i})^{m-1}
	  \pi(z_{i^{'}})\pi(z_{i^{'}})^{*}\\
        &&\> +(1-q^{2m})\sum_{k>i}\pi(z_{i})^{m-1}\pi(z_{k})\pi(z_{k})^{*};
\end{IEEEeqnarray*}
\item 
for $i+j < 2n+1$,  $i \neq j$,
\begin{displaymath}
\pi(z_i)^*\pi(z_j)^m = q^m\pi(z_j)^m\pi(z_i)^* + 
          mq^m(1-q^2)\epsilon_i\epsilon_jq^{\rho_i+
          \rho_j}\pi(z_j)^{m-1}\pi(z_{i^{'}})\pi(z_{j^{'}})^*;
\end{displaymath}
\item 
for $i> n$, 
\begin{displaymath}
\pi(z_{i})^*[\pi(z_{n}),\pi(z_{n+1})]^m= q^{2m}[\pi(z_{n}),\pi(z_{n+1})]^m \pi(z_{i})^*;
\end{displaymath}
\item
for $i > n$,
\begin{displaymath}
\pi(z_i)^*u_{\alpha_{2},\cdots ,\alpha_{2n-1},\alpha_0} = C u_{\alpha_{2},\cdots \alpha_{i-1},\alpha_i-1,\alpha_{i+1}\cdots ,\alpha_{2n-1},\alpha_0}
\end{displaymath}
where $C$ is some non-zero constant;
\item 
for $ n < i <2n$,
\begin{displaymath}
\pi(z_i)^*u_{\alpha_{2},\cdots \alpha_{i-1},0,\alpha_{i+1}\cdots ,\alpha_{2n-1},\alpha_0} = 0;
\end{displaymath}
\item
we have 
\begin{IEEEeqnarray*}{lCl}
\pi(z_{n})^*[\pi(z_{n}),\pi(z_{n+1})]^m
       &= &  q^{2m}[\pi(z_{n}),\pi(z_{n+1})]^m \pi(z_{n})^*+ (1-q^4)(1-q^2)\times\\
\IEEEeqnarraymulticol{3}{l}{\left(
      \sum_{l=0}^{k-1}q^{4l}[\pi(z_{n}),\pi(z_{n+1})]^{k-1-l}
      \pi(z_{n+1})[\pi(z_{n}),\pi(z_{n+1})]^l\right)\sum_{k>n+1}\pi(z_k)\pi(z_k)^*};
\end{IEEEeqnarray*}
\item
for $1<i \leq n$,
\begin{IEEEeqnarray*}{rCl}
\IEEEeqnarraymulticol{3}{l}{\pi(z_i)^*\pi(z_{i-1})^{\alpha_{i-1}}\cdots \pi(z_2)^{\alpha_2}[\pi(z_{n}),\pi(z_{n+1})]^{\alpha_0}u } \\
       &=&   C\pi(z_{i-1})^{\alpha_{i-1}}\cdots \pi(z_2)^{\alpha_2}\pi(z_i)^*[\pi(z_{n}),\pi(z_{n+1})]^{\alpha_0}u,
\end{IEEEeqnarray*}
where $C$ is some non-zero constant;
\item
we have
\begin{displaymath}
[\pi(z_{n}),\pi(z_{n+1})]^*[\pi(z_{n}),\pi(z_{n+1})]^m u =C[\pi(z_{n}),\pi(z_{n+1})]^{m-1}u,
\end{displaymath}
where $C$ is some non-zero constant;
\item 
for $1 <i\leq n$,
\begin{displaymath}
\pi(z_i)\pi(z{_i^{'}}^{m})=q^{2m}\pi(z_{i^{'}}^{m})\pi(z_i) -\sum_{k>i}(1-q^{2m})q^{i-k}\pi(z_k)\pi(z_{k^{'}})\pi(z_i^{m-1});
\end{displaymath}
\item
for $1 \leq i< n$,
\begin{displaymath}
\pi(z_i^*)[\pi(z_{n}),\pi(z_{n+1})]^m u =C\pi(z_i^{'})[\pi(z_{n}),\pi(z_{n+1})]^{m-1}u,
\end{displaymath}
where $C$ is some constant;
\item
we have
\begin{displaymath}
\pi(z_1)^*u_{\alpha_{2},\cdots ,\alpha_{n},0,\cdots,0,\alpha_0} = C u_{\alpha_{2},\cdots \alpha_{n},0,\cdots,0,\alpha_0-1}.
\end{displaymath}
\end{enumerate}
\eppsn 
\prf
We will prove part $4$ and part $9$ of this proposition. Other parts will follow by direct calculation using the commutation relations. 
\begin{enumerate}
\item
For $i> n+1$, it follows from $(\ref{c1})$. For $i=n+1$, it is enough to show for $m=1$.
\begin{IEEEeqnarray*}{lCl}
\IEEEeqnarraymulticol{3}{l}{\pi(z_{n+1})^*[\pi(z_{n}),\pi(z_{n+1})] }\\
        &=& \pi(z_{n+1})^*\pi(z_n)\pi(z_{n+1})-\pi(z_{n+1})^*\pi(z_{n+1})\pi(z_{n})\\
         &=& q^2\pi(z_n)\pi(z_{n+1})^*\pi(z_{n+1}) -
		  \pi(z_{n+1})\pi(z_{n+1})^*\pi(z_{n})\\
         &&\> -(1-q^2)\sum_{k>n+1}\pi(z_k)\pi(z_k^*)\pi(z_{n})\\
          &=&q^2\pi(z_n)\pi(z_{n+1})\pi(z_{n+1})^*+
                     q^2(1-q^2)\sum_{k>n+1}\pi(z_{n})\pi(z_k)\pi(z_k^*)\\
          &&\> -q^2\pi(z_{n+1})\pi(z_{n})\pi(z_{n+1})^*-q^2(1-q^2)
                         \sum_{k>n+1}\pi(z_{n})\pi(z_k)\pi(z_k^*)\\
           &=& q^2[\pi(z_{n}),\pi(z_{n+1})]\pi(z_{n+1})^*.
\end{IEEEeqnarray*}
\item
\begin{IEEEeqnarray*}{lCl}
\IEEEeqnarraymulticol{3}{l}{[\pi(z_{n}),\pi(z_{n+1})]^*[\pi(z_{n}),\pi(z_{n+1})]^{m}u } \\
         &=&  \pi(z_n^*)\pi(z_{n+1})^*[\pi(z_{n}),\pi(z_{n+1})]^{m}u +
		  \pi(z_{n+1}^*)\pi(z_n^*)[\pi(z_{n}),\pi(z_{n+1})]^{m}u \\
          &=&  C\pi(z_{n+1}^*)\sum_{l=0}^{m-1}q^{4l}[\pi(z_{n}),\pi(z_{n+1})]^{m-1-l}
		  \pi(z_{n+1})[\pi(z_{n}),\pi(z_{n+1})]^l u, \\
           && \mbox{(by part $7$ of  Proposition \ref{chap5-ppsn-2})}\\
          &=&  C\sum_{l=0}^{m-1}q^{4l}[\pi(z_{n}),\pi(z_{n+1})]^{m-1-l}
		  \pi(z_{n+1}^*)\pi(z_{n+1})[\pi(z_{n}),\pi(z_{n+1})]^l u\\
            && \mbox{(by   (\ref{c1}))}\\
            &=& C\sum_{l=0}^{m-1}q^{4l}[\pi(z_{n}),\pi(z_{n+1})]^{m-1-l}
			(\pi(z_{n+1})\pi(z_{n+1})^{*}  \\ 
           &&\> + (1-q^{2})\sum_{k>i}\pi(z_{k})\pi(z_{k})^{*})
		   [\pi(z_{n}),\pi(z_{n+1})]^l u\\
           &=&  C[\pi(z_{n}),\pi(z_{n+1})]^{m-1}u. 
\end{IEEEeqnarray*}
\end{enumerate} 
\qed

From  part $9$ of Proposition \ref{chap5-ppsn-2}, it follows that  $[\pi(z_{n}), \pi(z_{n+1})]^{\alpha_0}u \neq 0$.
Further we have  $\ker\pi(z_{i}) \subset \ker\pi(z_{2n}^*)=\left\{0\right\}$; hence 
 $u_{\alpha_{2},\cdots ,\alpha_{2n-1},\alpha_0}\neq 0$ for all 
 $(\alpha_{2},\cdots \alpha_{2n-1},\alpha_0) \in \bbn^{2n-1}$. Therefore  we can define
\begin{IEEEeqnarray}{rCl} \label{chap5-eqn-e}
e_{\alpha_{2},\cdots ,\alpha_{2n-1},\alpha_0} = \frac {u_{\alpha_{2},\cdots ,\alpha_{2n-1},\alpha_0}}
{\left\|u_{\alpha_{2},\cdots ,\alpha_{2n-1},\alpha_0}\right\|}.
\end{IEEEeqnarray}

\bppsn \label{chap5-ppsn-3}
Assume $\left\{e_{\alpha_{2},\cdots ,\alpha_{2n-1},\alpha_0}:(\sum_{l=2}^{2n-1}\alpha_{l})+2\alpha_0 \leq L\right\}$ form an 
orthonormal basis for $\clh_{\leq L} = \clh_{0}\oplus \clh_{1}\oplus \cdots \oplus \clh_{L}$.
If $2(r+s) +1\leq L$, then
\[
\left[\pi(z_{n}), \pi(z_{n+1})\right]^{r}\pi(z_{n+1})\left[\pi(z_{n}), \pi(z_{n+1})\right]^{s}u =
C\pi(z_{n+1})\left[\pi(z_{n}), \pi(z_{n+1})\right]^{r+s}u,
\]
where $C$ is a non-zero constant.
\eppsn 

\prf
It is enough to prove the statement for $r=1$.
The condition ensures that,
\[
 \left[\pi(z_{n}), \pi(z_{n+1})\right]\pi(z_{n+1})\left[\pi(z_{n}), \pi(z_{n+1})\right]^{s}u \in \clh_{\leq L}.
\]
 Hence
\begin{IEEEeqnarray*}{rCl}
\IEEEeqnarraymulticol{3}{l}{\left[\pi(z_{n}),  \pi(z_{n+1})\right]
               \pi(z_{n+1})\left[\pi(z_{n}), \pi(z_{n+1})\right]^{s}u \qquad\qquad} \\
      &=& \sum_{(\alpha_2\cdots \alpha_{2n-1},\alpha_0):
	          (\sum_{i=2}^{2n-1}\alpha_i)+  2\alpha_0 \leq L}
               C(\alpha_2, \cdots ,\alpha_0)e_{\alpha_2,\cdots ,\alpha_0},
\end{IEEEeqnarray*}
where 
\[
C(\alpha_2, \cdots ,\alpha_{2n-1},\alpha_0) = \left\langle \left[\pi(z_{n}),
\pi(z_{n+1})\right]\pi(z_{n+1})\left[\pi(z_{n}), \pi(z_{n+1})\right]^{s}u, e_{\alpha_2,\cdots ,\alpha_{2n-1},\alpha_0} \right\rangle.
\]
We will show that $C(\alpha_2,\cdots ,\alpha_{2n-1},\alpha_0) = 0 \mbox{ if } \alpha_{n+1}=1 \mbox{ and } \alpha_0 = s+1$. \\
\textbf{Case $1:$} $\alpha_{i} \neq 0$ for some $i > n+1.$\\ Applying part $4$ of the Proposition \ref{chap5-ppsn-2}, we get 
\[ 
\pi(z_i)^*\left[\pi(z_{n}), \pi(z_{n+1})\right]\pi(z_{n+1})\left[\pi(z_{n}), \pi(z_{n+1})\right]^{s}u = 0.
\] 
This shows that if $\alpha_{i} \neq 0$ for any $i \in \left\{n+2,n+3,\cdots ,2n-1\right\}$ then 
\begin{displaymath}
 C(\alpha_2, \cdots ,\alpha_{2n-1},\alpha_0) = 0.
 \end{displaymath}
\textbf{Case 2:} $\alpha_{n+1} \geq 1$  and  $\alpha_{i}= 0$  for all  $i > n+1$.
\begin{IEEEeqnarray*}{lCl}
     \IEEEeqnarraymulticol{3}{l}{\pi(z_{n+1})^*\left[\pi(z_{n}), \pi(z_{n+1})\right]
                 \pi(z_{n+1})\left[\pi(z_{n}), \pi(z_{n+1})\right]^{s}u}\\
      &=& q^2\left[\pi(z_{n}), \pi(z_{n+1})\right]\pi(z_{n+1}^*z_{n+1})
	             \left[\pi(z_{n}), \pi(z_{n+1})\right]^{s}u.\\
       &=& q^2\left[\pi(z_{n}), \pi(z_{n+1})\right]\pi(z_{n+1}z_{n+1}^*)
	   \left[\pi(z_{n}), \pi(z_{n+1})\right]^{s}u\\
        &&\> +\sum_{k> n+1}q^2(1-q^2)\left[\pi(z_{n}), \pi(z_{n+1})\right]
		\pi(z_kz_k^*)\left[\pi(z_{n}), \pi(z_{n+1})\right]^{s}u\\
         &=&q^{4s+2}(1-q^2)\left[\pi(z_{n}), \pi(z_{n+1})\right]^{s+1}u\\
		 &&\qquad (\mbox{since } u \in \bigcap_{i=1}^{2n-1} \mbox{ ker }\pi(z_{i}^{*})).
\end{IEEEeqnarray*} 
Now,
\begin{IEEEeqnarray}{lCl}
\left\langle u_{\alpha_{2},\cdots ,\alpha_{n+1},0,\cdots ,0,\alpha_0} ,\left[\pi(z_{n}), 
\pi(z_{n+1})\right]\pi(z_{n+1})\left[\pi(z_{n}), \pi(z_{n+1})\right]^{s}u \right\rangle \nonumber\\
=\left\langle u_{\alpha_{2},\cdots ,\alpha_{n+1}-1,0,\cdots ,0,\alpha_0} ,
\pi(z_{n+1}^*)\left[\pi(z_{n}), \pi(z_{n+1})\right]\pi(z_{n+1}\left[\pi(z_{n}), \pi(z_{n+1})\right]^{s}u \right\rangle.\nonumber\\
=\left\langle u_{\alpha_{2},\cdots ,\alpha_{n+1}-1,0,
\cdots ,0,\alpha_0} , q^{4s+2}(1-q^2)\left[\pi(z_{n}), \pi(z_{n+1})\right]^{s+1}u \right\rangle.\nonumber\\
 \begin{cases}
 \neq 0 & \mbox{ if } \alpha_{n+1} =1, \alpha_0= s+1,\alpha_{n-1}=\cdots \alpha_1=0, \cr
  = 0   & \mbox{ otherwise }. \cr
	\end{cases} \nonumber
\end{IEEEeqnarray} 
\textbf{Case 3:} $\alpha_{i}= 0 \mbox{ for all } i \geq n+1$.\\ 
By using the commutation relations, we have
\begin{IEEEeqnarray}{lCl}
\pi(z_n)^*u_{\alpha_{2},\cdots ,\alpha_{n},0,\cdots ,0,\alpha_0}\nonumber\\
= \pi(z_n)^{\alpha_{n}}\pi(z_n)^*\pi(z_{n-1})^{\alpha_{n-1}}\cdots, \left[\pi(z_{n}), \pi(z_{n+1})\right]^{\alpha_0}u \nonumber\\
\qquad + (1-q^{2\alpha_{n}})\sum_{k>n}\pi(z_n)^{\alpha_{n}-1}\pi(z_k)\pi(z_k)^*\pi(z_{n-1})^{\alpha_{n-1}}\cdots, \left[\pi(z_{n}),
\pi(z_{n+1})\right]^{\alpha_0}u.\nonumber\\
=C\pi(z_n)^{\alpha_{n}}\pi(z_{n-1})^{\alpha_{n-1}}\cdots, \pi(z_2)^{\alpha_2}\pi(z_n)^*\left[\pi(z_{n}), 
\pi(z_{n+1})\right]^{\alpha_0}u  \nonumber\\
\qquad +(1-q^{2\alpha_{n}})\sum_{k>n}\pi(z_n)^{\alpha_{n}-1}\pi(z_k)\pi(z_k)^*\pi(z_{n-1})^{\alpha_{n-1}}\cdots, 
\left[\pi(z_{n}), \pi(z_{n+1})\right]^{\alpha_0}u.\nonumber
\end{IEEEeqnarray} 
for some non-zero constant $C$. 

From part $5$ of the Proposition \ref{chap5-ppsn-2}, it follows that the
second term of the right hand side is 
$Cu_{\alpha_2,\cdots ,\alpha_{n-1},\alpha_{n}-1,0,\cdots ,0,\alpha_0} $.
Also, the first term of the right hand side is
\begin{multline*}
 C\pi(z_n)^{\alpha_{n}}\cdots, \pi(z_2)^{\alpha_2}[\pi(z_{n}),\pi(z_{n+1})]^{\alpha_0} \pi(z_{n})^*\\
 {}+C\pi(z_n)^{\alpha_{n}}\cdots, \pi(z_2)^{\alpha_2} (\sum_{l=0}^{\alpha-1}q^{4l}[\pi(z_{n}),\pi(z_{n+1})]^{\alpha_0-1-l}\pi(z_{n+1})[\pi(z_{n}),\pi(z_{n+1})]^l)u.
\end{multline*} 
Hence,
\begin{IEEEeqnarray}{lCl}
\pi(z_{n+1})^*\pi(z_n)^*u_{\alpha_{2},\cdots ,\alpha_{n},0,\cdots ,0,\alpha_0}\nonumber\\
 = C\pi(z_n)^{\alpha_{n}}\cdots \pi(z_2)^{\alpha_2}\pi(z_{n+1}^*)
 (\sum_{l=0}^{\alpha_0-1}q^{4l}[\pi(z_{n}),\pi(z_{n+1})]^{\alpha_0-1-l}
 \pi(z_{n+1})[\pi(z_{n}),\pi(z_{n+1})]^l)u\nonumber\\
=C\pi(z_n)^{\alpha_{n}}\cdots \pi(z_2)^{\alpha_2}
(\sum_{l=0}^{\alpha_0-1}q^{4l}
[\pi(z_{n}),\pi(z_{n+1})]^{\alpha_0-1-l}\pi(z_{n+1}^*)\pi(z_{n+1})[\pi(z_{n}),\pi(z_{n+1})]^l)u\nonumber\\
=C\pi(z_n)^{\alpha_{n}}\cdots \pi(z_2)^{\alpha_2}
(\sum_{l=0}^{\alpha_0-1}q^{4l}[\pi(z_{n}),\pi(z_{n+1})]^{\alpha_0-1-l}(\pi(z_{n+1})
\pi(z_{n+1}^*)\nonumber\\
 \quad +\sum_{k>n+1}\pi(z_k)\pi(z_k^*))[\pi(z_{n}),\pi(z_{n+1})]^l)u\nonumber\\
=C\pi(z_n)^{\alpha_{n}}\cdots \pi(z_2)^{\alpha_2}[\pi(z_{n}),\pi(z_{n+1})]^{\alpha_0-1}u.\nonumber
\end{IEEEeqnarray}
By the above calculation and by Proposition $\ref{chap5-ppsn-2}$, we have 
\begin{IEEEeqnarray}{lCl}
\left\langle u_{\alpha_2,\cdots ,\alpha_n,0.\cdots ,0,\alpha_0} , 
[\pi(z_{n}),\pi(z_{n+1})]\pi(z_{n+1})[\pi(z_{n}),\pi(z_{n+1})]^{s}u \right\rangle \nonumber\\
= \left\langle [\pi(z_{n}),\pi(z_{n+1})]^*u_{\alpha_2,\cdots ,\alpha_n,0.\cdots ,0,\alpha_0} , 
\pi(z_{n+1})[\pi(z_{n}),\pi(z_{n+1})]^{s}u \right\rangle\nonumber\\
= \left\langle \pi(z_{n+1}^*z_n^*-z_n^*z_{n+1}^*)u_{\alpha_2,\cdots ,\alpha_n,0.\cdots ,0,\alpha_0} ,
\pi(z_{n+1})[\pi(z_{n}),\pi(z_{n+1})]^{s}u \right\rangle\nonumber\\
= \left\langle \pi(z_{n+1}^*z_n^*)u_{\alpha_2,\cdots ,\alpha_n,0.\cdots ,0,\alpha_0} ,
\pi(z_{n+1})[\pi(z_{n}),\pi(z_{n+1})]^{s}u \right\rangle\nonumber\\
= \left\langle C\pi(z_n)^{\alpha_{n}}\pi(z_{n-1})^{\alpha_{n-1}}\cdots,
\pi(z_2)^{\alpha_2}[\pi(z_{n}),\pi(z_{n+1})]^{\alpha_0-1}u , \pi(z_{n+1})[\pi(z_{n}),\pi(z_{n+1})]^{s}u \right\rangle\nonumber\\
=\left\langle C\pi(z_{n+1}^*)\pi(z_n)^{\alpha_{n}}\pi(z_{n-1})^{\alpha_{n-1}}\cdots, 
\pi(z_2)^{\alpha_2}[\pi(z_{n}),\pi(z_{n+1})]^{\alpha_0-1}u , [\pi(z_{n}),\pi(z_{n+1})]^{s}u \right\rangle\nonumber\\
=0 \nonumber
\end{IEEEeqnarray}
This proves the claim.\qed

\blmma \label{chap5-lmma-ONB}
 Let $\pi$ be an irreducible representation on a Hilbert space $\clh$ with $\pi(z_{2n}) \neq 0$. Then
 $\left\{e_{\alpha_{2},\alpha_{2},\cdots ,\alpha_{2n-1},\alpha_0} : (\alpha_{2},\alpha_{3},
 \cdots \alpha_{2n-1},\alpha_0) \in \bbn^{2n-1}\right\}$ defined in equation $(\ref{chap5-eqn-e})$, form an orthonormal basis for $\clh$.
\elmma
\prf 
From Corollary \ref{cr2}, it is enough to show that for $\alpha \neq \beta$, $u_{\alpha}$ is orthogonal to $u_{\beta}$. 
We apply induction on $L_{\alpha}:=(\sum_{i=2}^{2n-1}\alpha_{i}) + 2\alpha_0$. For $L_{\alpha} = 0$, 
claim is true as $u \neq 0$. Assume the hypothesis for $L_{\alpha} \leq N-1$. Note that $u_{\alpha_{2},\cdots ,\alpha_{2n-1},\alpha_0}
\in \clh_{L_{\alpha}}$. Hence, by induction hypothesis and Corollary $\ref{cr2}$, it follows that 
$\left\{e_{\alpha_{2},\cdots ,\alpha_{2n-1},\alpha_0}: (\sum_{l=2}^{2n-1}\alpha_{l})+2\alpha_0=m \right\}$ form 
an orthonormal basis of $\clh_{m}$ for $m \leq N-1$.

If $\alpha$ and $\beta$ are such that $L_{\alpha} \neq L_{\beta},$ then $u_{\alpha} \in \clh_{L_{\alpha}}$ and 
$u_{\beta} \notin \clh_{L_{\alpha}}$, which shows that $u_{\alpha}$ and $u_{\beta}$ are orthogonal. Take $\alpha$ 
and $\beta$ such that $L_{\alpha}= L_{\beta} = N$. Assume  $\alpha_{i} \neq 0$ for some $i >n$. Choose maximum such $i$.
From part $6$ of Proposition \ref{chap5-ppsn-2}, it follows that
\begin{IEEEeqnarray}{rCl}
\Big\langle u_{\alpha} , u_{\beta} \Big\rangle & = & 
\Big\langle u_{\alpha_2,\cdots ,\alpha_{i-1},\alpha_i-1,0,\cdots ,0,\alpha_0} , \pi(z_i^*)u_{\beta}\Big\rangle \nonumber\\
 & = & \Big\langle u_{\alpha_2,\cdots ,\alpha_{i-1},\alpha_i-1,0,\cdots ,0,\alpha_0} ,
 Cu_{\beta_2,\cdots ,\beta_{i-1},\beta_{i}-1,\beta_{i+1},\cdots , \beta_{2n-1},\beta_0} \Big\rangle \nonumber
\end{IEEEeqnarray}
where $C$ is a non zero constant. Now, by using induction we get
$\left\langle u_{\alpha} , u_{\beta} \right\rangle \neq 0$ if and only if  
$\alpha = \beta$.
Hence it is enough to consider $\alpha$ and $\beta$ such that 
$\alpha_{i} = \beta_{i} =0$ for $i>n$. Let $\alpha_{n} \neq 0$.
\begin{IEEEeqnarray*}{rCl}
\IEEEeqnarraymulticol{3}{l}{\pi(z_n^*)u_{\beta_2,\cdots ,\beta_{n},0,\cdots ,0,\beta_0}}\\
     &=& (\pi(z_n^{\beta_n})\pi(z_n^*)+ 
                 \sum_{k>i}C\pi(z_n^{\beta_n-1})\pi(z_k)\pi(z_k^*))
                 u_{\beta_2,\cdots ,\beta_{n-1},0,\cdots ,0,\beta_0} \\
     &=&C \pi(z_n^{\beta_n})\pi(z_{n+1})
	    u_{\beta_2,\cdots ,\beta_{n-1},0,\cdots ,0,\beta_0-1}+
               Cu_{\beta_2,\cdots ,\beta_{n}-1,0,\cdots ,0,\beta_0}  \\
       & & \mbox{ (by Proposition~\ref{chap5-ppsn-3}) }   \\
     &=& \left(C\pi(z_{n+1})\pi(z_n^{\beta_n})+\sum_{k>n+1}C\pi(z_k)
	 \pi(z_k^{'})\pi(z_n^{\beta_n-1})\right)
         u_{\beta_2,\cdots ,\beta_{n-1},0,\cdots ,0,\beta_0-1} \\
& &\quad +\> Cu_{\beta_2,\cdots ,\beta_{n}-1,0,\cdots ,0,\beta_0} .
\end{IEEEeqnarray*}
Hence
\begin{IEEEeqnarray*}{rCl}
\IEEEeqnarraymulticol{3}{l}{\left\langle u_{\alpha} , u_{\beta} \right\rangle }\\
& = & \left\langle u_{\alpha_2,\cdots ,\alpha_{n-1},\alpha_{n}-1,0,\cdots ,0,\alpha_0} , \pi(z_n^*)u_{\beta} \right\rangle \\
& = &  \left\langle u_{\alpha_2,\cdots ,\alpha_{n-1},\alpha_{n}-1,0,\cdots ,0,\alpha_0} , 
Cu_{\beta_{2},\cdots ,\beta_{n}-1,0,\cdots ,0,\beta_0}
+ C\pi(z_{2n})\pi(z_1)u_{\beta_{2},\cdots ,\beta_{n}-1,0,\cdots ,0,\beta_0-1} \right\rangle \\
& = & \left\langle u_{\alpha_2,\cdots ,\alpha_{n-1},\alpha_{n}-1,0,\cdots ,0,\alpha_0} 
, Cu_{\beta_{2},\cdots ,\beta_{n}-1,0,\cdots ,0,\beta_0}\right\rangle \nonumber\\
& & \quad +  \left\langle \pi(z_1^*)u_{\alpha_2,\cdots ,\alpha_{n-1},\alpha_{n}-1,0,\cdots ,0,\alpha_0} ,
Cu_{\beta_{2},\cdots ,\beta_{n}-1,0,\cdots ,0,\beta_0-1}\right\rangle \nonumber\\
&=& \left\langle u_{\alpha_2,\cdots ,\alpha_{n-1},\alpha_{n}-1,0,\cdots ,0,\alpha_0} , 
Cu_{\beta_{2},\cdots ,\beta_{n}-1,0,\cdots ,0,\beta_0}\right\rangle \nonumber\\
& & \quad +\left\langle u_{\alpha_2,\cdots ,\alpha_{n-1},\alpha_{n}-1,0,\cdots ,0,\alpha_0-1} ,
Cu_{\beta_{2},\cdots ,\beta_{n}-1,0,\cdots ,0,\beta_0-1}\right\rangle.
\end{IEEEeqnarray*}
Again induction proves the claim. So, we will consider $\alpha$ and $\beta$ such that $\alpha_{i} = \beta_{i} =0$ for $i \geq n$.
Assume that for some $i \in \left\{2,3,\cdots,n-1\right\}, \alpha_{i} \neq 0$ or $\beta_{i} \neq 0$. 
Choose maximum such $i$. Without loss of generality, we assume that $\alpha_i \neq 0$. 
\begin{IEEEeqnarray*}{lCl}
\IEEEeqnarraymulticol{3}{l}{\pi(z_{i}^*)
                  u_{\beta_{2},\cdots ,\beta_{i},0,\cdots ,0,\beta_0}}\\
  &=& (\pi(z_{i})^{\beta_{i}}\pi(z_{i})^{*} +q^{2\rho_{i}}(1-q^{2\beta_{i}})\pi(z_{i})^{\beta_{i}-1}
       \pi(z_{i^{'}})\pi(z_{i^{'}})^{*} \\
   &&+\>(1-q^{2\beta_{i}})\sum_{k>i}\pi(z_{i})^{\beta_{i}-1}
            \pi(z_{k})\pi(z_{k})^{*})\pi(z_{i-1})^{\beta_{i-1}}\cdots
			\pi(z_2)^{\beta_2}[\pi(z_{n}),\pi(z_{n+1})]^{\beta_0}u\\
  &=& C\pi(z_{i})^{\beta_{i}}\pi(z_{i}^{'})\pi(z_{i-1})^{\beta_{i-1}}
                     \cdots \pi(z_2)^{\beta_2}[\pi(z_{n}),\pi(z_{n+1})]^{\beta_0-1}u
                +Cu_{\beta_{2},\cdots ,\beta_{i}-1,0,\cdots ,0,\beta_0} \\
   &&+\sum_{i<k\leq n}C\pi(z_{i})^{\beta_{i}-1}\pi(z_{k})\pi(z_{k})^{*})
            \pi(z_{i-1})^{\beta_{i-1}}\cdots ,\pi(z_2)^{\beta_2}
               [\pi(z_{n}),\pi(z_{n+1})]^{\beta_0-1}u.\\
&=&  C\pi(z_{i}^{'})\pi(z_{i})^{\beta_{i}}\pi(z_{i-1})^{\beta_{i-1}}\cdots \pi(z_2)^{\beta_2}[\pi(z_{n}),\pi(z_{n+1})]^{\beta_0-1}u 
+Cu_{\beta_{2},\cdots ,\beta_{i}-1,0,\cdots ,0,\beta_0}\\
&& +\sum_{k > i^{'}}C\pi(z_{k}))\pi(z_{k^{'}})\pi(z_{i})^{\beta_{i}-1}
             \pi(z_{i-1})^{\beta_{i-1}}
\cdots ,\pi(z_2)^{\beta_2}[\pi(z_{n}),\pi(z_{n+1})]^{\beta_0-1}u\\
&&   +\sum_{i<k \leq n}C\pi(z_{i})^{\beta_{i}-1}\pi(z_{k^{'}}))\pi(z_{k})
             \pi(z_{i-1})^{\beta_{i-1}}
             \cdots ,\pi(z_2)^{\beta_2}[\pi(z_{n}),\pi(z_{n+1})]^{\beta_0-1}u\\
&=&Cu_{\beta_{2},\cdots ,\beta_{i}-1,0,\cdots ,0,\beta_0} \\
&& +\sum_{n <k \leq 2n }C\pi(z_{k}))\pi(z_{k^{'}})\pi(z_{i})^{\beta_{i}-1}
         \pi(z_{i-1})^{\beta_{i-1}}
\cdots ,\pi(z_2)^{\beta_2}[\pi(z_{n}),\pi(z_{n+1})]^{\beta_0-1}u.
\end{IEEEeqnarray*}
Hence
\begin{IEEEeqnarray*}{lCl}
\IEEEeqnarraymulticol{3}{l}{\left\langle u_{\alpha}, u_{\beta} \right\rangle }\\
&=&  \left\langle u_{\alpha_2,\cdots ,\alpha_{i-1},\alpha_{i}-1,0,\cdots ,0,\alpha_0}, 
   \pi(z_{i}^*)u_{\beta_{2},\cdots ,\beta_{i},0,\cdots ,0,\beta_0}\right\rangle\\
& = &  \left\langle u_{\alpha_2,\cdots ,\alpha_{i-1},\alpha_{i}-1,0,\cdots ,0,\alpha_0} ,
      Cu_{\beta_{2},\cdots ,\beta_{i}-1,0,\cdots ,0,\beta_0} 
            + C\pi(z_{2n})\pi(z_1)
       u_{\beta_{2},\cdots ,\beta_{i}-1,0,\cdots ,0,\beta_0-1} \right\rangle \\
& = & \left\langle u_{\alpha_2,\cdots ,\alpha_{i-1},\alpha_{i}-1,0,\cdots ,0,\alpha_0} , 
Cu_{\beta_{2},\cdots ,\beta_{i}-1,0,\cdots ,0,\beta_0}\right\rangle \\
& &  \quad +\left\langle \pi(z_1^{*})u_{\alpha_2,\cdots ,\alpha_{i-1},\alpha_{i}-1,0,\cdots ,0,\alpha_0} ,
Cu_{\beta_{2},\cdots ,\beta_{i}-1,0,\cdots ,0,\beta_0-1}\right\rangle\\
&=& \left\langle u_{\alpha_2,\cdots ,\alpha_{i-1},\alpha_{i}-1,0,\cdots ,0,\alpha_0} ,
Cu_{\beta_{2},\cdots ,\beta_{i}-1,0,\cdots ,0,\beta_0}\right\rangle \\
& & \quad +\left\langle u_{\alpha_2,\cdots ,\alpha_{n-1},\alpha_{i}-1,0,\cdots ,0,\alpha_0-1} , 
Cu_{\beta_{2},\cdots ,\beta_{i}-1,0,\cdots ,0,\beta_0-1}\right\rangle.
\end{IEEEeqnarray*}
Again induction will settle the claim. Now we take $\alpha$ and $\beta$ such that $\alpha_{i} = \beta_{i} =0$ for all  $i \neq 0$. Then from part $9$ of  
Proposition~\ref{chap5-ppsn-2}, it follows that 
\[
[\pi(z_{n}),\pi(z_{n+1})]^*[\pi(z_{n}),\pi(z_{n+1})]^{\beta_{0}}u 
=  C [\pi(z_{n}),\pi(z_{n+1})]^{\beta_0-1}. 
\]
Hence
\begin{IEEEeqnarray}{lCl}
\left\langle u_{\alpha}, u_{\beta} \right\rangle 
& = & \left\langle [\pi(z_{n}),\pi(z_{n+1})]^{\alpha_{0}}u , [\pi(z_{n}),\pi(z_{n+1})]^{\beta_{0}}u \right\rangle \nonumber \\
& = & \left\langle [\pi(z_{n}),\pi(z_{n+1})]^{\alpha_{0}-1}u , [\pi(z_{n}),\pi(z_{n+1})]^*[\pi(z_{n}),\pi(z_{n+1})]^{\beta_{0}}u \right\rangle \nonumber\\ 
& = & \left\langle [\pi(z_{n}),\pi(z_{n+1})]^{\alpha_{0}-1}u , C[\pi(z_{n}),\pi(z_{n+1})]^{\beta_{0}-1}u \right\rangle \nonumber
\end{IEEEeqnarray} 
This completes the proof. \qed
\bcrlre \label{cr1}
 If $\pi$ and $\pi^{'}$ are two irreducible representations of $C(H_{q}^{2n})$ with $\pi(z_{2n}) \neq 0 \mbox{ and } \pi^{'}(z_{2n}) \neq 0$, then
  for all $\alpha \in \bbn^{2n-1}$, we have
\[ \left\|u_{\alpha}\right\| = \left\|u_{\alpha}^{'}\right\|.
\]
where $u_{\alpha}$ and $u_{\alpha}^{'}$ are defined as above.
\ecrlre 
Now we aim to find all irreducible representations of $C(H_{q}^{2n})$. One way is to do explicit calculations
to determine the operators $z_{1},z_{2},\cdots ,z_{2n}$ as done in case of the odd dimensional quantum spheres. 
But in this case, computations are more complicated. So, to avoid complicated calculations,
we show that one can completely determine an irreducible representation  $\pi$ of $C(H_{q}^{2n})$ given that
$\pi(\omega) \neq 0$ and $\pi(z_{2n})u = tu$ for some fixed  $t \in \bbbt$. Then we use representation of the quantum Stiefel manifold 
$C(SP_{q}(2n)/SP_{q}(2n-2))$ to get explicit description of the representation.

\bthm
Let $\pi$ and $\pi^{'}$ be  irreducible representations of $C(H_{q}^{2n})$ on a Hilbert space $\clh$ and $\clh^{'}$ respectively 
such that $\pi(z_{2n})|_{\bigcap_{i=1}^{2n-1}\ker \pi(z_i^*)} = tI = \pi^{'}(z_{2n})|_{\bigcap_{i=1}^{2n-1}\ker \pi^{'}(z_i^*)} $ for
$t \in \bbbt$. Then $\pi$ and $\pi^{'}$ are equivalent.
\ethm 
\prf
Without loss of generality, we can take $t=1$. Let $u$ and $u^{'}$ are unit vectors 
in $\bigcap_{i=1}^{2n-1}\ker \pi(z_i^*)$ and $\bigcap_{i=1}^{2n-1}\ker \pi^{'}(z_i^*)$ respectively. 
From  Lemma \ref{chap5-lmma-ONB}, we have  canonical orthonormal bases for $\clh$ and $\clh^{'}$ given by 
\[
\left\{e_{\alpha_{2},\alpha_{3},\cdots ,\alpha_{2n-1},\alpha_0}:  (\alpha_{2},\alpha_{3},\cdots \alpha_{2n-1},\alpha_0) \in \bbn^{2n-1}\right\}
\]
and 
\[
\left\{e_{\alpha_{2},\alpha_{3},\cdots ,\alpha_{2n-1},\alpha_0}^{'}: 
(\alpha_{2},\alpha_{3},\cdots \alpha_{2n-1},\alpha_0) \in \bbn^{2n-1}\right\}
\]
 respectively. Define $U:\clh \rightarrow \clh^{'}$ by
\begin{displaymath}
U: e_{\alpha_{2},\alpha_{3},\cdots ,\alpha_{2n-1},\alpha_0}\longmapsto e_{\alpha_{2},\alpha_{3},\cdots ,\alpha_{2n-1},\alpha_0}^{'}
\end{displaymath}
From Corollary \ref{cr1}, $U(u_{\alpha_{2},\alpha_{3},\cdots ,\alpha_{2n-1},\alpha_0})= u_{\alpha_{2},\alpha_{3},\cdots ,\alpha_{2n-1},\alpha_0}^{'}$. 
We know, $\clh = \oplus_{m \in \bbn} \clh_m$ and $\clh^{'} = \oplus_{m \in \bbn} \clh_m^{'}$.
where $\clh_{m}$ and $\clh_{m}^{'}$ are the eigenspaces of $\pi(\omega)$ and $\pi^{'}(\omega)$ respectively,
corresponding to the eigenvalue $q^{2m}$. Clearly $U(\clh_{m}) = \clh_{m}^{'}$.  
We need to show that, $U\pi(z_{i})U^{*}=\pi^{'}(z_{i})$, or equivalently $U\pi(z_{i}^{*})U^{*}=\pi^{'}(z_{i}^{*})$ for all 
$i \in \left\{1,2,\cdots ,2n\right\}$. We split the proof into several parts.\\
$\textbf{(A)}\quad$  $i=2n$:
\begin{displaymath}
\pi(z_{2n}^*)u_{\alpha_{2},\cdots ,\alpha_{2n-1},\alpha_0}=q^{(\sum_{l=2}^{2n-1}\alpha_{l})+2\alpha_{0}}
u_{\alpha_{2},\cdots ,\alpha_{2n-1},\alpha_0}.
\end{displaymath}
\begin{displaymath}
\pi^{'}(z_{2n}^*)u_{\alpha_{2},\cdots ,\alpha_{2n-1},\alpha_0}=q^{(\sum_{l=2}^{2n-1}\alpha_{l})+2\alpha_{0}}
u_{\alpha_{2},\cdots ,\alpha_{2n-1},\alpha_0}^{'}.
\end{displaymath}
$\textbf{(B)}\quad$  $n<i<2n$:
\begin{displaymath}
\pi(z_{i}^*)u_{\alpha_{2},\cdots ,\alpha_{2n-1},\alpha_0} = Cu_{\alpha_{2},\cdots ,\alpha_{i-1},\alpha_{i},\cdots ,\alpha_{2n-1},\alpha_0}
\end{displaymath}
\begin{displaymath}
\pi^{'}(z_{i}^*)u_{\alpha_{2},\cdots ,\alpha_{2n-1},\alpha_0}^{'} = 
Cu_{\alpha_{2},\cdots ,\alpha_{i-1},\alpha_{i},\cdots ,\alpha_{2n-1},\alpha_0}^{'}.
\end{displaymath}
Note that the constant $C$ is same in both equations. Hence we have,
\begin{displaymath}
U\pi(z_{i}^*)U^{*}=\pi^{'}(z_{i}) \qquad \mbox{ for } n< i \leq 2n.
\end{displaymath}
$\textbf{(C)}\quad$  $i= n$:\\
In this case,
we will use induction on the dimension of the eigenspaces of $\pi(\omega)$. 
For $m=0$, we have $\pi(z_{i}^{*})u = 0= \pi{'}(z_{i}^{*})u^{'}$.  
Assume that $U\pi(z_{i})U^{*}_{|\clh_{\leq m}} =\pi^{'}(z_{i})_|\clh_{\leq m}^{'}$.
Take $ u_{\alpha_{2},\cdots ,\alpha_{2n-1},\alpha_0} \in \clh_{m+1}$.\\
\textbf{Case $1:$} $\alpha_{j} \neq 0$ for some  $j > n$,  and  $\alpha_{k} =0$ for all $k>j$. 
\bean
\pi(z_{n}^{*})u_{\alpha_{2},\cdots ,\alpha_{j},0,\cdots ,0,\alpha_0} & = & C\pi(z_j)\pi(z_n^{*})
u_{\alpha_{2},\cdots ,\alpha_{j}-1,0,\cdots ,0,\alpha_0}, \\
&&  \qquad  \mbox{(by  (\ref{c1})}\\
& = &  C \pi(z_j)U^*\pi^{'}(z_n^{*})u_{\alpha_{2},\cdots ,\alpha_{j}-1,0,\cdots ,0,\alpha_0}^{'},\\
&& \qquad  \mbox{ (by induction) }\\
& = &  C U^*\pi{'}(z_j)UU^*\pi^{'}(z_n^{*})u_{\alpha_{2},\cdots ,\alpha_{j}-1,0,\cdots ,0,\alpha_0}^{'}, \\
& = &  C U^*\pi{'}(z_j)\pi^{'}(z_n^{*})u_{\alpha_{2},\cdots ,\alpha_{j}-1,0,\cdots ,0,\alpha_0}^{'}, \\
& = &  \pi{'}(z_{n})u_{\alpha_{2},\cdots ,\alpha_{j},0,\cdots ,0,\alpha_0}^{'}.
\eean 
\textbf{Case $2:$} $\alpha_{j} = 0$ for all $j>n$.\\ 
From part $7$, part $8$ of the Proposition \ref{chap5-ppsn-2} and Proposition \ref{chap5-ppsn-3}, we have
\bean
\pi(z_{n}^*)u_{\alpha_{2},\cdots ,\alpha_{n},0,\cdots ,0,\alpha_0} & = & Cu_{\alpha_{2},\cdots ,\alpha_{n}-1,0,\cdots ,0,\alpha_0} 
+ Cu_{\alpha_{2},\cdots ,\alpha_{n},1,0,\cdots ,0,\alpha_0-1}.\\
\pi^{'}(z_{n}^*)u_{\alpha_{1},\cdots ,\alpha_{n},0,\cdots ,0,\alpha_0}^{'} 
& = & Cu_{(\alpha_{2},\cdots ,\alpha_{n}-1,0,\cdots ,0,\alpha_0}^{'} + Cu_{(\alpha_{2},\cdots ,\alpha_{n},1,0,\cdots ,0,\alpha_0-1)}^{'}.\\
\eean
Hence we get
 \begin{displaymath}
U\pi(z_n)^*U^* = \pi^{'}(z_n).
\end{displaymath}
$\textbf{(D)}\quad$  $1 < i < n$: \\
\textbf{Case $1:$} $\alpha_{j} \neq 0$ for some  $j > i$,  and  $\alpha_{k} =0$ for all $k>j$. \\
This follows exactly as in the case $i = n$. \\
\textbf{Case $2:$} $\alpha_{j} = 0$, for all $j \geq i$. \\
It follows from part $2$ 
of the Proposition \ref{chap5-ppsn-2} and by using the fact $U\pi(z_k)^*U^* = \pi^{'}(z_k)$ for all $k>i$.\\
\textbf{Case $3:$}  $\alpha_{j} = 0$, for all $j \geq i$.\\
 From part $8$ and part $11$ of the Proposition \ref{chap5-ppsn-2}, we have
\bean 
\pi(z_{i})^*u_{\alpha_{2},\cdots ,\alpha_{i-1},0,\cdots ,0,\alpha_0} 
& = & Cu_{\alpha_{2},\cdots ,\alpha_{i-1},0,\cdots ,0,\underbrace{1}_{i^{'}-1 \mbox{ th place }},0,\cdots,0,\alpha_0}\\
\pi^{'}(z_{i})^*u_{\alpha_{2},\cdots ,\alpha_{i-1},0,\cdots ,0,\alpha_0}
& = &  Cu_{\alpha_{2},\cdots ,\alpha_{i-1},0,\cdots ,0,\underbrace{1}_{i^{'}-1 \mbox{ th place }},0,\cdots,0,\alpha_0}^{'}.\\
\eean
which settles the claim for $1< i < n$.\\
$\textbf{(E)}\quad$ For $i=1$:\\
We again use induction.
Take $u_{\alpha}$ such that  $\alpha_{j} \neq 0$ for some $j \neq 0$. Choose $j$ to be max $\{k: \alpha_k \neq 0\}$.
\bean
\pi(z_1)^*u_{\alpha_2\cdots ,\alpha_{j},0,\cdots ,0,\alpha_0} & = & C\pi(z_j)\pi(z_1)^*u_{\alpha_2,\cdots ,\alpha_{j}-1,0,\cdots ,0,\alpha_0},\\
 & = & CU^*\pi{'}(z_j)UU^*\pi{'}(z_1)^*u_{\alpha_2,\cdots ,\alpha_{j}-1,0,\cdots ,0,\alpha_0}^{'},\\
& = & CU^*\pi{'}(z_j)\pi{'}(z_1)^*u_{\alpha_2,\cdots ,\alpha_{j}-1,0,\cdots ,0,\alpha_0}^{'},\\
& = & \pi^{'}(z_1)^*u_{\alpha_2\cdots ,\alpha_{j},0,\cdots ,0,\alpha_0}^{'}.\\
\eean 
For $\alpha$ such that $\alpha_j = 0$, for all $j \neq 0$, it follows from part $12$ of the Proposition \ref{chap5-ppsn-2} and induction. 
Hence, we have $U\pi(z_i^*)U^* = \pi^{'}(z_i^*)$ for  $1 \leq i \leq 2n$, which proves the claim. \qed

We will now discuss the general case. Let $\pi$ be an irreducible representation of $C(H_{q}^{2n})$ on a Hilbert space $\clh$ such 
that $\pi(z_{2n})=\pi(z_{2n-1})=\cdots =\pi(z_{k+1})=0, \mbox{ and } \pi(z_k) \neq 0$. It follows from (\ref{c8}) that $z_k$ is normal. 
Denote $z_k^*z_k$ by $\omega$. By the same reasoning,  $\clh$ decomposes as
\[\clh = \oplus_{m \in \bbn}\clh_{m}. 
\] 
where $\clh_{m}$ is the eigenspace of $\pi(\omega)$ corresponding to the eigenvalue $q^{2m}$.

Let $K$ be a subspace of $\bigcap_{i=1}^{k-1} \mbox{ ker }\pi(z_{i}^{*})$ 
invariant under the $C^*$-algebra generated by $\pi(z_{k})$. Define
\begin{eqnarray*}
\clh^{K} = \mbox{ linear span }\Big\{\pi(z_{1})^{\alpha_{1}}\pi(z_{2})^{\alpha_{2}}\cdots \pi(z_{k-1})^{\alpha_{k-1}}h: h \in K\Big\}.
\end{eqnarray*}
In the same way, one can show that, $\clh^{K}$ is an invariant subspace of a representation $\pi$ and  hence 
by irreducibility of the representation, $\bigcap_{i=1}^{k-1} \mbox{ ker }\pi(z_{i}^{*})$ is one dimensional. 
Pick $u \in \bigcap_{i=1}^{k-1} \mbox{ ker }\pi(z_{i}^{*}) = \clh_{0}$. As $\pi(z_k)$ keeps $\clh_{0}$ invariant 
and  $\pi(z_k)_{|\clh_0}$ is a unitary operator, we get $\pi(z_k)u= tu \mbox{ for some } t \in \bbbt$.
\bthm
Let $1 \leq k \leq 2n$ and let $\pi$ be an irreducible representation of $C(H_{q}^{2n})$ on a Hilbert space $\clh$ 
such that $\pi(z_{2n})=\pi(z_{2n-1})=\cdots =\pi(z_{k+1})=0$. Assume $\pi(z_{k})u = tu, t \in \bbbt$, where $u$
is defined as above. Then $\pi$ is the unique representation up to equivalence which satisfies these conditions.
\ethm   
\prf
This is essentially the previous proof with some minor modifications.\\
\textbf{Case $1:k>n.$} \\
Define
\begin{IEEEeqnarray*}{lCl}
\IEEEeqnarraymulticol{3}{l}{u_{\alpha_{1},\alpha_{2},\cdots ,\alpha_{2n-k},\alpha_{2n-k+2},\cdots,\alpha_{k-1},\alpha_0}}\\
&=& \pi(z_1)^{\alpha_1}\cdots \pi(z_{2n-k})^{\alpha_{2n-k}}  
   \pi(z_{k-1})^{\alpha_{k-1}}\cdots \pi(z_{2n-k+2})^{\alpha_{2n-k+2}}[\pi(z_{n}),
\pi(z_{n+1})]^{\alpha_{0}}.
\end{IEEEeqnarray*}  
where $\alpha_{i} \in \bbn$. (Note the differences between the definition  of $u_{\alpha}$ given above and that given  in equation (\ref{chap5-eqn-u})).
\[
e_{\alpha_{1},\alpha_{2},\cdots ,\alpha_{2n-k},\alpha_{2n-k+2},\cdots,\alpha_{k-1},\alpha_0} 
= \frac {u_{\alpha_{1},\alpha_{2},\cdots ,\alpha_{2n-k},\alpha_{2n-k+2},\cdots,\alpha_{k-1},\alpha_0}}
{\left\|u_{\alpha_{1},\alpha_{2},\cdots ,\alpha_{2n-k},\alpha_{2n-k+2},\cdots,\alpha_{k-1},\alpha_0}\right\|}\\
\]   
Similar calculations as have been done in Lemma \ref{chap5-lmma-ONB}  will prove that 
\[
\left\{e_{\alpha_{1},\cdots ,\alpha_{2n-k},\alpha_{2n-k+2},\cdots,\alpha_{k-1},\alpha_0}, : 
(\alpha_{1},\cdots ,\alpha_{2n-k},\alpha_{2n-k+2},\cdots,\alpha_{k-1},\alpha_0) \in \bbn^{k-1}\right\}
\]
 form an orthonormal basis for $\clh$. By the same argument as used in the previous theorem, one  proves that 
the representation satisfying  $\pi(z_{2n})=\pi(z_{2n-1})=\cdots =\pi(z_{k+1})=0$ and $\pi(z_{k})u = tu$ is unique.\\
\textbf{Case $2:k\leq n$}. \\
First observe that the relations satisfied by $\pi(z_1),\pi(z_2),\cdots ,\pi(z_k)$ are same 
as the defining relations of the odd dimensional quantum sphere $S_q^{2k+1}$ for which we know that the claim holds.
(Note that we can also proceed as in the previous case  and establish the claim.) \qed

We have shown so far that if there exists an irreducible representation $\pi$ such that
$\pi(z_{2n})=\pi(z_{2n-1})=\cdots =\pi(z_{k+1})=0, \mbox{ and } \pi(z_{k})u = tu \mbox{ for } t \in T$,
then it is unique. Existence of these representations still needs to be shown.

By Theorem \ref{quotient}, $C(SP_{q}(2n)/SP_{q}(2n-2))$ is the $C^{*}$-subalgebra of $C(SP_{q}(2n))$ generated by
$\left\{u_{m}^{2n}\right\}_{m=1}^{2n}$. Now, if we look at the relations $I_{st}^{ij}$ involving $u_{m}^{2n}$ and $u_{m}^{1}$ 
by putting  $(i,j)=(1,1)$ and $(2n,1)$, we get the relations satisfied by generators of $C(H_q^{2n})$ where $z_{m} = u_{2n+1-m}^{2n}$.
From the universal property of $C(H_{q}^{2n})$, there exists a map $\eta :C(H_{q}^{2n}) \longrightarrow C(SP_{q}(2n)/SP_{q}(2n-2))$ such that 
$\eta(z_{j}) = u_{2n+1-j}^{2n}$ for all $j \in \left\{1,2,\cdots ,2n\right\}$.

Denote by $\omega_{k}$ the following word in the Weyl group of $sp_{2n}$,
\[
\omega_k =\begin{cases}
             I & \mbox{ if } k =1, \cr
             s_1s_2\cdots s_{k-1} & \mbox{ if } 2 \leq k \leq n,\cr
             s_1s_2\cdots s_{n-1}s_{n}s_{n-1}\cdots s_{2n-k+1} & \mbox{ if } n < k \leq 2n.\cr
						\end{cases}
\] 
Let $\eta_{t,\omega_k} = \pi_{t,\omega_{k}}\circ \eta$. Hence we have an irreducible representation 
$\eta_{t,\omega_k}$ of $C(H_{q}^{2n})$ such that 
$\eta_{t,\omega_k}(z_{2n})=\eta_{t,\omega_k}(z_{2n-1})=\cdots =\eta_{t,\omega_k}(z_{k+1})=0$ and $\eta_{t,\omega_k}(z_{k})u = tu$ 
where $ 1< k \leq 2n$. This gives an explicit description of the irreducible representations satisfying these conditions. 

For $k=1$, define $\eta_{t,I}:C(H_{q}^{2n})\rightarrow \bbc$ by  $\eta_{t,I}(z_j) = t\delta_{1j}$.
The set $\left\{\eta_{t,I}: t \in T\right\}$ gives all one dimensional irreducible representations of $C(H_{q}^{2n})$. 
Also, it satisfies $\eta_{t,I}(z_{2n})=\eta_{t,I}(z_{2n-1})=\cdots =\eta_{t,I}(z_2)=0 \mbox{ and } \eta_{t,I}(z_1)u = tu$.   
\bcrlre \label{repn}
The set $\left\{\eta_{t,\omega_k}: 1\leq k \leq 2n, t \in \bbbt \right\}$ gives a complete list of irreducible representations of $C(H_{q}^{2n})$.
\ecrlre 
To get a faithful representation of $C(H_{q}^{2n})$, define $\eta_{\omega_k}: C(H_{q}^{2n})\rightarrow C(\bbbt) \otimes \scrt^{\otimes k-1}$ 
by  $\eta_{\omega_k}(a)(t) = \eta_{t,\omega_k}(a) \quad$ for all $a \in  C(H_{q}^{2n})$.
\bcrlre
$\eta_{\omega_{2n}}$ is a faithful representation of $C(H_{q}^{2n})$.
\ecrlre
\prf It is easy to see that any irreducible representation factors through $\eta_{\omega_{2n}}$ as, $\omega_k$ is a subword of $\omega_{2n}$.
This proves the claim.
\qed
\bthm 
The homomorphism $\eta : C(H_{q}^{2n})\rightarrow C(SP_{q}(2n)/SP_{q}(2n-2))$ is an isomorphism.              
\ethm
\prf Clearly $\eta$ is a surjective homomorphism as $\eta(z_i)$ are generators of $C(SP_{q}(2n)/SP_{q}(2n-2))$. 
It follows from Corollary \ref{repn} that all irreducible representations of $C(H_{q}^{2n})$  factor through $\eta$ which
shows that $\eta$ is injective. This proves the claim.
\qed

We separate out some important facts that will be useful in determining the $K$- groups of $C(H_{q}^{2n})$.
\bcrlre \label{cr5}
Let $C_{1}= C(\bbbt)$ and for $ 2 \leq k \leq 2n$, let $C_{k} = \eta_{\omega_{k}}(C(H_{q}^{2n}))$.
Then the set $\left\{\eta_{t,\omega_l}: 1\leq l \leq k, t \in \bbbt \right\}$ gives a complete list of irreducible representations of $C_k$. 
\ecrlre

\bcrlre \label{cr6}
 Let $\pi = \eta_{t,\omega_k}$ . Then
\begin{enumerate}
\item 
for $1 \leq k \leq n$, one has
\begin{IEEEeqnarray}{lCl}
\pi(z_1)^{\alpha_1}\cdots \pi(z_{k-1})^{\alpha_{k-1}}1_{\left\{1\right\}}(\pi(z_{k}^{*}z_{k}))= 
Cp_{\alpha_{1},0} \otimes p_{\alpha_{2},0} \otimes ...\otimes p_{\alpha_{k-1},0}.\nonumber
\end{IEEEeqnarray}
\item
for $n< k \leq 2n$, one has
\begin{IEEEeqnarray*}{lCl}
\IEEEeqnarraymulticol{3}{l}{\pi(z_1)^{\alpha_1}\cdots \pi(z_{2n-k})^{\alpha_{2n-k}}
 \pi(z_{k-1})^{\alpha_{k-1}}\pi(z_{k-2})^{\alpha_{k-2}}
    \cdots \pi(z_{2n-k+2})^{\alpha_{2n-k+2}}}\\
 \IEEEeqnarraymulticol{3}{l}{[\pi(z_{n}),\pi(z_{n+1})]^{\alpha_{0}}
 1_{\left\{1\right\}}(\pi(z_{k}^{*}z_{k}))}\\
&=& Cp_{\alpha_{1},0} \otimes \cdots p_{\alpha_{2n-k},0}\otimes p_{\alpha_{2n-k+2},0}\otimes \cdots p_{\alpha_{n},0} \otimes p_{\alpha_{0},0} 
\otimes p_{\alpha_{n+1},0} \otimes \cdots \otimes p_{\alpha_{k-1},0},
\end{IEEEeqnarray*}
 where $p_{i,j}$ is the rank one operator on $\ell^2{2}(\bbn)$ sending the basis element $e_{j}$ to $e_{i}$ and $C$ is some non-zero constant.
\end{enumerate}
\ecrlre

\newsection{$K$-groups of $C(H_{q}^{2n})$}

  Neshveyev \& Tuset (\cite{NesTus-2012ab}) proved the $KK$-equivalence of the $C^*$-algebra $C(G/K)$ and $C(G_q/K_q)$. 
  Hence one can determine the $K$-groups of $C(G_q/K_q)$ from those of $C(G/K)$ via the equivalence.
  As a consequence, generators of the $K$-groups of $C(G_q/K_q)$  would be images of generators of
  the corresponding $K$-groups of $C(G/K)$ under the equivalence. However the proof of the equivalence that is known 
  is  existential in nature and it is hard to find any generator using that equivalence. 
  Here we obtain the $K$-groups of $C(SP_{q}(2n)/SP_{q}(2n-2))$ in a computationally more tractable way 
  and give an explicit description of the generators of the $K$-groups.

 We first derive certain exact sequences analogous to that for the odd dimensional  quantum sphere (see~\cite{She-1997ac}). 
 We then apply the six-term sequence in $K$-theory to compute the $K$- groups of $C(SP_{q}(2n)/SP_{q}(2n-2))$. 
  Let $p_{i,j}$ be the rank one operator on $L_{2}(\bbn)$ sending  basis element 
 $e_{j}$ to $e_{i}$ and $p$ be the operator $p_{0,0}$. 
\blmma \label{lk1}
Let $C_{1}= C(\bbbt)$ and for $2 \leq k \leq 2n$, let  $C_{k} = \eta_{\omega_{k}}(C(H_{q}^{2n}))$.
Then $C(\bbbt) \otimes \clk(\ell^2(\bbn))^{\otimes (k-1)}$ is contained in $C_{k}$. 
Moreover, for $2 \leq k \leq 2n$ we have the exact sequence,
\[
      0 \longrightarrow C(\bbbt) \otimes \clk(\ell^2(\bbn))^{\otimes(k-1)} \longrightarrow 
      C_{k} \stackrel{\sigma_{k}}{\longrightarrow} C_{k-1} \longrightarrow 0.
\] 
where $\sigma_{k}$ is the restriction of $(1^{\otimes (k-1)}\otimes \sigma)$ to $C_{k}$ and $\sigma : \scrt \rightarrow \bbc$ 
is the homomorphism such that $\sigma(S)=1$.
\elmma
\prf 
First we prove that $C(\bbbt) \otimes \clk(\ell^2(\bbn))^{\otimes(k-1)}$ is contained in $C_{k}$. For $k \leq n$, 
it follows from a result of Sheu (\cite{She-1997ac}, Theorem $4$) as $C_{k}$ is isomorphic to $C(S_q^{2k-1})$. For $k > n$, and $m \geq 0$, 
\begin{displaymath}
\eta_{\omega_{k}}(z_{k}^{m}1_{\left\{1\right\}}(z_{k}^{*}z_{k})) = t^{m} \otimes
	\underbrace{p \otimes p\otimes ...\otimes p}_{k-1},\quad 
  \eta_{\omega_{k}}(z_{k}^{*m}1_{\left\{1\right\}}(z_{k}^{*}z_{k})) = t^{-m} \otimes \underbrace{p \otimes p\otimes ...\otimes p}_{k-1}.
\end{displaymath}
Also from Corollary $\ref{cr6}$, it follows that 
\begin{IEEEeqnarray*}{lCl}
\IEEEeqnarraymulticol{3}{l}{\eta_{\omega_{k}}\bigl((z_1)^{m_1}\cdots (z_{2n-k})^{m_{2n-k}}(z_{k-1})^{m_{k-1}}(z_{k-2})^{m_{k-2}}
\cdots (z_{2n-k+2})^{m_{2n-k+2}}}\\
\IEEEeqnarraymulticol{3}{l}{[(z_{n}),(z_{n+1})]^{m_{0}}1_{\left\{1\right\}}(z_{k}^{*}z_{k})\bigr)} \\
&=&  Ct^{(\sum_{i=0,i \neq 2n-k+1}^{k-1}m_{i})+m_0} \\
&&p_{m_1,0} \otimes \cdots p_{m_{2n-k},0}\otimes 
p_{m_{2n-k+2}}\otimes \cdots p_{m_{n},0} \otimes p_{m_{0},0} \otimes p_{m_{n+1},0} \otimes \cdots \otimes p_{m_{k-1},0}.
\end{IEEEeqnarray*}
which shows that $t \otimes p_{m_{1},0} \otimes p_{m_{2},0} \otimes ...\otimes p_{m_{k-1},0}$  and 
$1 \otimes p_{m_{1},0} \otimes p_{m_{2},0} \otimes ...\otimes p_{m_{k-1},0}$ are contained in $C_{k}$. Hence
$C_{k}$ contains $C(\bbbt) \otimes \clk(\ell^2(\bbn))^{\otimes(k-1)}$.

	                     It is easy to see that  $\sigma_{k}$ vanishes on $C(\bbbt) \otimes \clk(\ell^{2}(\bbn))^{\otimes(k-1)}$. 
	                     Also, it follows from Corollary \ref{cr5} that any irreducible representation of 
	                     $C_{k}$ is of the form $\eta_{t,\omega_l}$ where $l \leq k$ and $t\in \bbbt$.
	                     Hence an irreducible representation of $C_{k}$ that  vanishes on
	                     $C(\bbbt) \otimes \clk(\ell^{2}(\bbn))^{\otimes(k-1)}$ is of the form
	                     $\eta_{t,\omega_l}$ where $l \leq k-1$ and $t\in \bbbt$. But this factors 
	                     through $\sigma_{k}$. Thus we get the desired exact sequence. \qed
											
\brmrk Neshveyev \& Tuset (\cite{NesTus-2012ab}) obtained a composition series for $C(G_q/K_q)$ for 
any Poisson-Lie closed subgroup $K$ of $G$. In particular, when $K$ is $C(SP_q(2n-2))$, one  gets a 
composition series for $C(SP_{q}(2n)/SP_{q}(2n-2))$. Note that the series of exact sequence derived in the Lemma \ref{lk1} is 
different from that given in \cite{NesTus-2012ab}.
\ermrk
Define, for $1 \leq k \leq 2n$, 
\begin{displaymath}
 u_{k} = t \otimes
	\underbrace{p \otimes p\otimes ...\otimes p}_{k-1} + 1 - 1 \otimes
	\underbrace{p \otimes p\otimes ...\otimes p}_{k-1}.
\end{displaymath}
  It is easy to check that 
	$u_{k} =\eta_{\omega_{k}}(z_{k}1_{\left\{1\right\}}(z_{k}^{*}z_{k}) +1- 1_{\left\{1\right\}}(z_{k}^{*}z_{k}))$. 
	Therefore $u_{k}$ is an unitary operator which is contained in 
	$C_{k}$.
\bthm
Let $1 \leq k \leq 2n$.The $K$-groups $K_{0}(C_{k})$ and $K_{1}(C_{k})$ are both isomorphic to $\bbz$ and in particular,
$[u_{k}]$ form a $\bbz$-basis for  $K_{1}(C_{k})$ and $[1]$ form a $\bbz$-basis for  $K_{0}(C_{k})$. 
\ethm
\prf 
We apply induction on $k$. For $k=1$, this is clear. Assume the result is true for $k-1$. From Lemma \ref{lk1}, 
we have the short exact sequence
\[
  0 \longrightarrow C(\bbbt) \otimes \clk(\ell^{2}(\bbn))^{\otimes(k-1)} \longrightarrow C_{k} 
  \stackrel{\sigma_{k}}{\rightarrow} C_{k-1} \longrightarrow 0.
\]
which gives rise to the following six-term sequence in $K$-theory. 

\begin{tikzpicture}[node distance=1cm,auto]
	\tikzset{myptr/.style={decoration={markings,mark=at position 1 with %
	    {\arrow[scale=2,>=stealth]{>}}},postaction={decorate}}}
\node (A){$ K_0(C(\bbbt) \otimes \clk(\ell^{2}(\bbn))^{\otimes (k-1)})$};
\node (B)[node distance=4cm, right of=A]{$K_0(C_{k})$};
\node (Up)[node distance=2cm, right of=B][label=above:$K_{0}(\sigma_{k})$]{};
\node (C)[node distance=4cm, right of=B]{$K_0(C_{k-1})$};
\node (D)[node distance=2cm, below of=C]{$K_1(C(T) \otimes \clk(\ell^{2}(\bbn))^{\otimes (k-1)})$};
\node (E)[node distance=4cm, left of=D]{$K_1(C_{k})$};
\node (F)[node distance=4cm, left of=E]{$K_1(C_{k-1})$};
\draw[myptr](A) to (B);
\draw[myptr](B) to (C);
\draw[myptr](C) to node{{ $\delta$}}(D);
\draw[myptr](D) to (E);
\draw[myptr](E) to (F);
\draw[myptr](F) to node{{ $\partial$}}(A);
\end{tikzpicture}

To compute the six term sequence, we determine $\delta$ and $\partial$. Since $\sigma_{k}(1) = 1$, it follows that $\delta([1]) = 0$. 
Also, the operator $\widetilde{X} = t \otimes \underbrace{q^{N} \otimes q^{N} \otimes...}_{k-2} \otimes S^{*}$ is in $C_{k}$ as 
$\eta_{\omega_{k}}(z_{k-1}) - \widetilde{X}$ lies in $C(\bbbt) \otimes \clk(\ell^{2}(\bbn))^{\otimes(k-1)}$. Let
\[
X = 1_{\left\{1\right\}}(\widetilde{X}^{*}\widetilde{X})\widetilde{X} +1-1_{\left\{1\right\}}(\widetilde{X}^{*}\widetilde{X}).
\]	
Then $X$ is an isometry such that $\sigma_{k}(X) = u_{k-1}$ and hence
\[
\partial ([u_{k-1}]) = [1-X^{*}X]-[1-XX^{*}]= [1 \otimes \underbrace{p \otimes p\otimes ...\otimes p}_{k-1}].
\]
Now by the Kunneth theorem for the tensor product of $C^{*}$-algebras (see \cite{Bla-1998aa} ), 
it follows that $C(\bbbt) \otimes \clk(\ell^{2}(\bbn))^{\otimes(k-1)}$ has $K_{0}$-group isomorphic to $\bbz$ generated by $[1 \otimes
	\underbrace{p \otimes p\otimes ...\otimes p}_{k-1}]$ and $K_{1}$-group  isomorphic to $\bbz$ generated by $[u_{k}]$.
	Induction hypothesis and the above calculation shows that $\partial$ is an isomorphism and hence $K_{0}(i)$ is the zero map.
	Therefore
	$K_{0}(\sigma_{k})$ is injective. Since $\delta$ is zero, $K_{0}(\sigma_{k})$ is surjective. Hence 
	$K_{0}(C_{k})$ is isomorphic to $\bbz$ and is generated by $[1]$. Similarly $K_{1}(i)$ is injective as $\delta$ is the zero map. 
	Also, since $\partial$ is an isomorphism, $K_{1}(\sigma_{k})$ is the zero map. This shows that $K_{1}(i)$ is surjective. 
	Hence $K_{1}(C_{k})$ is isomorphic to $\bbz$ and is generated by $[u_{k}]$.
	This establishes the claim. \qed

\noindent\begin{footnotesize}\textbf{Acknowledgement}:
I would like to thank Prof. Arup Kumar Pal, my supervisor, for his constant support. I would also like to thank S. Sundar 
for useful discussions on various topics.
\end{footnotesize}


\noindent{\sc Bipul Saurabh} (\texttt{bipul9r@isid.ac.in})\\
         {\footnotesize Indian Statistical
Institute, 7, SJSS Marg, New Delhi--110\,016, INDIA}

\end{document}